\documentclass[final,3p,onecolumn,times,numbers,sort&compress]{elsarticle}

\usepackage{amssymb}
\usepackage[fleqn]{amsmath}
\usepackage{algorithm}
\usepackage{algpseudocode}
\usepackage{graphicx}
\usepackage{booktabs}
\usepackage{siunitx}
\usepackage{hyperref}
\usepackage{float}
\usepackage{enumitem}
\setlist{leftmargin=15pt,labelindent=15pt}
\setlist[enumerate]{wide=0pt, leftmargin=15pt, labelwidth=15pt, align=left}

\hypersetup{
  colorlinks   = true, 
  urlcolor     = blue, 
  linkcolor    = blue, 
  citecolor   = blue 
}

\journal{Elsevier}

\begin{document}

\begin{frontmatter}



\title{Topology optimization for particle flow problems using Eulerian-Eulerian model with a finite difference method}


\author[a]{Chih-Hsiang Chen}
\author[a]{Kentaro Yaji\corref{c1}}
\cortext[c1]{Corresponding author}
\ead{yaji@mech.eng.osaka-u.ac.jp}
\affiliation[a]{organization={Department of Mechanical Engineering, Graduate School of Engineering, Osaka University},
            addressline={2-1, Yamadaoka}, 
            city={Suita},
	   state={Osaka},
            postcode={ 565-0871}, 
            country={Japan}}

\begin{abstract}
Particle flow processing is widely employed across various industrial applications and technologies. Due to the complex interactions between particles and fluids, designing effective devices for particle flow processing is challenging. In this study, we propose a topology optimization method to design flow fields that effectively enhance the resistance encountered by particles. Particle flow is simulated using an Eulerian-Eulerian model based on a finite difference method. Automatic differentiation is implemented to compute sensitivities using a checkpointing algorithm. We formulate the optimization problem as maximizing the variation of drag force on particles while reducing fluid power dissipation. Initially, we validate the finite difference flow solver through numerical examples of particle flow problems and confirm that the corresponding topology optimization produces a result comparable to the benchmark problem. Furthermore, we investigate the effects of Reynolds and Stokes numbers on the optimized flow field. The numerical results indicate that serpentine flow fields can effectively enhance the variation in particle drag force.
\end{abstract}


\end{frontmatter}


\section{Introduction}

The process of handling particle flow plays a pivotal role in a wide range of industrial applications~\cite{klaewkla2011review, duan2015green}. Among the factors governing these processes, the particle residence time---the duration particles spend in a specific region---serves as a fundamental factor of process efficiency, affecting the degree of reaction~\cite{duduta2011semi, wang2020review}, mixing quality~\cite{santana2019development, suh2010review}, and heat transfer effect~\cite{chaurasia2020thermal, das2006heat} in industrial devices. A fundamental approach to adjusting particle residence time is to modify the resistance acting on the particles. This can be achieved by controlling operational conditions~\cite{hong2017millifluidic, kehlenbeck2002particle}, altering particle properties~\cite{phirommark2023cfd, alhamdan1997residence, lee2018parametric}, or designing flow paths to impose different levels of resistance~\cite{afrin2014modeling, al2020study, bovskovic2011residence}. Furthermore, it is beneficial to increase the resistance encountered by particles for applications that require extending particle residence time. For instance, in fluidized beds, baffles are used to control particle movement by reducing back-mixing and redistributing particles, which increases resistance and enhances mixing efficiency~\cite{yang2015experimental, kreesaeng2024effect, rennebaum2024effect}. In microreactors, complex channel designs are employed to extend flow paths and ensure complete chemical reactions~\cite{an2012computational, coyle2008micro, peres2019analysis}. Additionally, porous media are incorporated within particle heating receivers to increase the resistance to particle flow, allowing the particles to absorb solar energy effectively and thereby improving thermal efficiency~\cite{ho2016review, ho2019sun}. However, the complex interactions between particles and fluid make it challenging to design an effective flow path to increase particle flow resistance. Although size optimization is commonly employed for such requirements~\cite{aubin2009effect}, predefined shapes offer limited flexibility in determining optimized configurations at the conceptual design stage. 

Unlike size and shape optimization, topology optimization can generate optimized configurations without the designer's intuition at the conceptual design stage. The homogenization design method for topology optimization, first proposed by Bendsøe and Kikuchi, processes heterogeneous materials into homogeneous forms to achieve optimized material distribution for enhanced structural performance~\cite{bendsoe1988generating}. Moreover, the density approach, also known as the SIMP (Solid Isotropic Material with Penalization) method, employs a continuous function to enable a smooth transition of design variables throughout the entire design domain in the optimization process~\cite{bendsoe2008topology}. The use of a penalization factor in the SIMP method eliminates intermediate values of design variables, improving material distribution from a manufacturing viewpoint.

Based on the SIMP method, topology optimization for Stokes flow was pioneered by Borrvall and Petersson~\cite{borrvall2003topology}. They introduced a fictitious force term to model porous domains, allowing for no-slip boundary conditions by setting a sufficiently high value of fictitious force. Subsequent research has broadened the scope of topology optimization in fluid dynamics. For instance, Guest et al. developed a methodology for optimizing creeping flow conditions~\cite{guest2006topology}. Gersborg-Hansen et al. proposed a topology optimization method for laminar channel flow problems utilizing Navier-Stokes equations~\cite{gersborg2005topology}. Topology optimization methods for unsteady flow problems have been extended by Kreissl et al.~\cite{kreissl2011topology} and Deng et al.~\cite{deng2011topology}. Owing to the high degree of design freedom, research in fluid topology optimization is not only limited to pure fluid issues but also includes a wide range of studies that consider interactions with various physical fields~\cite{yaji2015topology, yaji2018topology,chen2019computational, jenkins2015level, li2023topology}.

In the case of multi-physics problems, there has been growing interest in the topology optimization of particle flows in recent years. Andreasen proposed a topology optimization approach for optimizing inertial microfluidic devices for particle manipulation~\cite{andreasen2020framework}. Yoon proposed studies focused on controlling particle trajectories in fluid~\cite{yoon2022transient,yoon2021development,choi2023matlab}. These studies employ the Eulerian-Lagrangian model to explore the topology optimization of single or multiple particles in fluids. It is crucial to highlight that particle modeling can be conducted using either the Eulerian-Lagrangian approach or the Eulerian-Eulerian approach, each offering distinct advantages for particle flow analysis. In the Eulerian-Lagrangian approach, the fluid is treated as a continuum, and each particle is tracked individually within the Lagrangian framework, enabling the detailed simulation of particle trajectories and interactions with the fluid phase~\cite{zahari2018introduction}. Conversely, the Eulerian-Eulerian approach treats both the gas and particles as continuous phases with their conservation equations, allowing for the investigation of behaviors of numerous particles, including sedimentation, fluidization, and particle residence time distribution~\cite{balakin2010eulerian,zhang2019simulation,gao2012review}. Although prior studies have utilized the Eulerian-Lagrangian approach to optimize particle manipulations and trajectory problems, employing the Eulerian-Lagrangian approach to simulate the collective behavior of particles can lead to high consumption of computational resources due to the necessity of tracking each particle individually~\cite{ariyaratne2016cfd}. Therefore, it becomes important to optimize problems that address issues related to the behavior of numerous particles using the Eulerian-Eulerian approach.
 
On the other hand, according to the literature review~\cite{alexandersen2020review}, most studies on topology optimization employ the finite element method, while a few studies utilize the finite volume method~\cite{yan2023topology}, lattice Boltzmann method~\cite{makhija2012topology} and particle-based method~\cite{sasaki2019topology}. Although the finite element method is well-suited for structural problems with arbitrary shapes of objects, it requires relatively significant memory resources~\cite{zawawi2018review} and involves specific treatments for conservation law in fluid problems~\cite{fanxi2017fast,barman2016introduction}. Despite its challenges with complex geometries, the finite difference method effectively handles fluid simulation~\cite{anderson1995computational,moukalled2016finite} and supports high-order schemes for complex fluid behavior~\cite{cheng2009high,rizzetta2008high}. From the view of standard topology optimization methods, an optimized structure is typically determined by the distribution of design variables. There is no need to employ unstructured meshes to address issues of arbitrary structural shapes, which makes the finite difference method an appropriate choice for fluid topology optimization.

In this paper, we present a topology optimization method for particle flow to increase the resistance acting on the particles effectively. For the numerical model, we employ a finite difference method and the Eulerian-Eulerian approach to simulate particle flow behaviors. The finite difference scheme called the NAPPLE (Nonstaggered Artificial Pressure for Pressure-Linked Equation) algorithm is selected as our flow solver. The characteristic of the NAPPLE algorithm lies in its use of specific assumptions that allow all state variables to be defined on a collocated grid, thereby simplifying programming and enhancing computational efficiency~\cite{lee1992artificial}. The optimization problems are solved via the GCMMA (Globally Convergent Method of Moving Asymptotes)~\cite{CCSA}. The sensitivity information is computed with automatic differentiation, and a checkpointing algorithm is utilized to reduce memory requirements~\cite{griewank2008evaluating,griewank2000algorithm}. In the presented numerical examples, we examine how variations in the Reynolds number and the Stokes number influence the configuration of the optimized flow field.

The subsequent sections of this paper are arranged as follows: In Section~\ref{Mathematical model}, we define the governing equations and boundary conditions. In Section~\ref{Problem formulation}, we detail the approach for topology optimization and present a problem designed to maximize changes in particle drag force. In Section~\ref{Numerical implementation}, we first describe the NAPPLE algorithm in detail, followed by the demonstration of the overall numerical implementation, which includes sensitivity calculation and the optimization algorithm. In Section~\ref{Results and discussion}, some numerical results are discussed to confirm the efficacy of this proposed method. Finally, the summarization and future work are provided in Section~\ref{Conclusion}.

\section{Mathematical model}
\label{Mathematical model}
\subsection{Assumption}
The interactions between each phase play an important role in particle flow behavior. In the condition of low particle volume fraction, the effect of particles on the fluid is negligible, a concept known as one-way coupling. With the increases in the particle volume fraction, the effect of two-way coupling should be considered for flow characteristics, which means the fluid and particles mutually influence each other. When the volume fraction of particles exceeds 0.1, the particle flow is regarded as dense particle flow, where the particle-particle interactions need to be taken into account~\cite{crowe2011multiphase}.

The assumption of one-way coupling is suitable for applications with dilute particles, such as microfluidic devices~\cite{kim2010new, wang2019industry, khashan2014computational} or small-scale particle heating receivers~\cite{hicdurmaz2022numerical}. Therefore, we adopt a one-way coupling condition with sufficient diluteness of the particles. In addition, due to the complex interactions between turbulence and particles, and given the widespread applications of particle flow in laminar regimes~\cite{sajeesh2014particle}, we focus on steady-state laminar flow in this study. Some basic assumptions are adopted for simplification, as follows:

\begin{itemize}
\item The fluid phase is treated as a steady incompressible laminar flow.
\item Due to the dilute condition of the particles, the influence of particles on the fluid and the interactions between particles are not considered.
\item The shape effect is not considered, and particles are spherical.
\item Since the lift force on particles has a minor effect in the investigated cases, the particle lift force is negligible~\cite{kotoky2018development}.
\end{itemize}

\subsection{Governing equations}
\label{Governing equations}
\begin{figure}[h]
\centering
\includegraphics[width=1.0\textwidth]{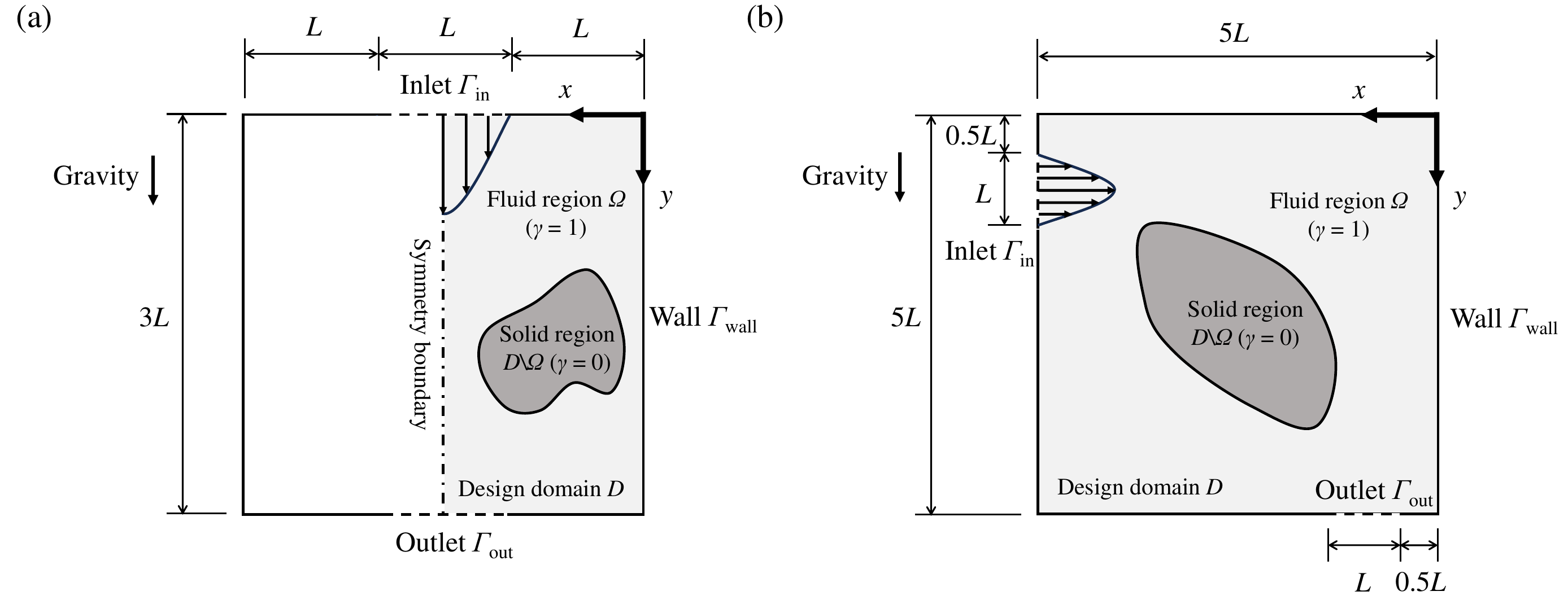} 
\caption{Schematic diagrams for the analysis domains and boundary conditions: (a) symmetrical flow (b) asymmetrical flow.}
\label{Fig.6}
\end{figure}

We adopt the Eulerian-Eulerian model as our core approach for the numerical investigation. Both phases are formulated by the conservation equations for continuity and momentum. 
The governing equations for the fluid phase are defined as follows:
\begin{align}\label{gef}
&\nabla \cdot ({\phi_{\rm f} \rho_{\rm f} \bold{u}_{\rm f} \bold{u}_{\rm f}}) =  -\phi_{\rm f} \nabla p + \nabla \cdot (\mu \phi_{\rm f}(\nabla \bold{u}_{\rm f} + \nabla \bold{u}^\intercal_{\rm f})) + \phi_{\rm f} \rho_{\rm f} \bold{g}, \\[10pt]
\label{fc}
&\nabla \cdot (\phi_{\rm f} \bold{u}_{\rm f}) = 0,
\end{align}
where $\phi_{\rm f}$ is the fluid volume fraction,  $\rho_{\rm f}$ is the fluid density, $\mu$ is the fluid viscosity, $\bold{g}$ is the gravity constant, $p$ is the pressure, and $\bold{u_{\rm f}}$ is the velocity vector of fluid phase. Similarly, the governing equations for the particle phase in the Eulerian-Eulerian formulation are defined as follows:
\begin{align}\label{gep}
&\nabla \cdot ({\phi_{\rm p} \rho_{\rm p} \bold{u}_{\rm p} \bold{u}_{\rm p}}) =  -\phi_{\rm p} \nabla p + \phi_{\rm p} \rho_{\rm p} \bold{g} + K(\bold{u_{\rm f}} - \bold{u_{\rm p}}), \\[10pt]
\label{pc}
&\nabla \cdot (\phi_{\rm p} \bold{u}_{\rm p}) = 0,
\end{align}
where $\bold{u_{\rm p}}$ is the velocity vector of particle phase, $\phi_{\rm p}$ is the particle volume fraction, $\rho_{\rm p}$ is the particle density, $K$ is the interphase momentum exchange coefficient. To consider the momentum exchange between the fluid phase and the particle phase,  the source term related to the momentum exchange coefficient $K$ is incorporated into the momentum equation for the particle phase, which can be determined by the Gidaspow model~\cite{gidaspow1994multiphase} expressed as follows:
\begin{equation}
K = 
\begin{cases}
\begin{aligned}
&\frac{3}{4} C_{\rm d}  \frac{\phi_{\rm p} \phi_{\rm f} \rho_{\rm f} \lvert \bold{u}_{\rm f} - \bold{u}_{\rm p}\rvert} {d_{\rm p}} \phi_{\rm f}^{-2.65} , \text{ if $\phi_{\rm f} >$ 0.8} \\
&150\frac{\phi_{\rm p}^2 \mu_{\rm f}}{\phi_{\rm f}d_{\rm p}^2} + 1.75  \frac{\phi_{\rm p} \rho_{\rm f} \lvert \bold{u}_{\rm f} - \bold{u}_{\rm p}\rvert} {d_{\rm p}} , \text{ if $\phi_{\rm f} \leq $ 0.8},
\end{aligned}
\end{cases}
\end{equation}
where $d_{\rm p}$ is the particle diameter, and $C_{\rm d}$ is the drag coefficient, which is defined by Kolve~\cite{gidaspow1994multiphase} as follows:
\begin{equation}
C_{\rm d} =
\begin{cases}
\begin{aligned}
&\dfrac{24}{Re_{\rm p}}  					,\quad \text{if $Re_{\rm p} < $ 1} 		\\
&\dfrac{24}{Re_{\rm p}}(1 + 0.15Re_{\rm p}^{0.687}) 	,\quad \text{if 1$ \leq Re_{\rm p} \leq $ 1000},
\end{aligned}
\end{cases}
\end{equation}
where $Re_{\rm p}$ is the particle Reynolds number, which can be described as:
\begin{equation}
Re_{\rm p}  = \frac{\rho_{\rm f} \lvert \bold{u}_{\rm f} - \bold{u}_{\rm p}\rvert d_{\rm p}}{\mu}.
\end{equation}
In this study, we consider the cases of symmetric and asymmetric flow conditions. Fig.~\ref{Fig.6} shows the analysis domains and boundary conditions. For momentum conservation, the Dirichlet boundary conditions for velocities of both phases are imposed at the inlet. The Dirichlet boundary condition for pressure is specified at the outlet, and no-slip boundary conditions are applied to the walls. The expressions are shown below:
\begin{alignat}{4}
\label{bc1}
&\bold{u}_{\rm f} = \bold{u}^{\rm in}_{\rm f}, \quad \bold{u}_{\rm p} = \bold{u}^{\rm in}_{\rm p} 			&\text{\quad  on \quad}&\mathit{\Gamma}_{\rm in}, 		\\
\label{bc12}
&\bold{u}_{\rm f} = \bold{0}, \quad \bold{u}_{\rm p} = \bold{0}			 					&\text{\quad  on \quad}&\mathit{\Gamma}_{\rm wall}, 		\\
\label{bcp}
&p = p^{\rm out}  															&\text{\quad  on \quad}&\mathit{\Gamma}_{\rm out},
\end{alignat}
where $\bold{u}^{\rm in}_{\rm f}$ and $\bold{u}^{\rm in}_{\rm p}$ are the fluid velocity and particle velocity at the inlet, respectively, and $p^{\rm out}$ is the prescribed value for pressure at the outlet. Besides, for the volume fraction, the Dirichlet boundary condition for the particle phase is also defined on the inlet, and the Neumann boundary condition is specified on the remaining boundaries. The expressions are as follows:
\begin{alignat}{4}
\label{bcp2}
&\bold{\phi}_{\rm p} = \phi^{\rm in}_{\rm p}, 		&\text{\quad  on \quad}&\mathit{\Gamma}_{\rm in}, 				\\
\label{bc2}
&\nabla {\phi}_{\rm p} \cdot \bold{n} = 0,			&\text{\quad  on \quad}&\mathit{\Gamma}_{\rm wall}, \mathit{\Gamma}_{\rm out}, 
\end{alignat}
where $\phi^{\rm in}_{\rm p}$ is the prescribed value of particle volume fraction at the inlet. Note that since the summation of the volume fraction of both phases is equal to one at any given node, specifying the particle volume fraction on the boundaries automatically determines the fluid volume fraction. Both numerical examples also account for gravity, aligned in the direction of the $y$-axis.
\section{Problem formulation}
\label{Problem formulation}
\subsection{Basic concept of topology optimization}
\label{Basic concept of topology optimization}
Topology optimization aims to derive an optimized design within a design domain $D$, guided by a specific objective function and constraints. Essentially, 
the basic idea of topology optimization involves introducing design variables to transform a structural design problem into a material distribution problem. We formulate the expression of design variables as follows:
\begin{equation}
\bold{\chi}(\bold{x}) = 
\begin{cases}
1 , \text{ if $\bold{x} \in \mathit{\Omega}$}\\
0, \text{ if $\bold{x} \in D \backslash \mathit{\Omega}$},
\end{cases}
\end{equation}
where $\bold{x}$ is the position in $D$, and $\mathit{\Omega}$ represents the flow domain in $D$.  With the discontinuity of the characteristic function
$\bold{\chi}(\bold{x})$, relaxation techniques are commonly employed to achieve a continuous transition for numerical treatment. In this study, a popular method called the density approach is utilized for topology optimization, where the material distribution is defined by a continuous function with values ranging from $0 \leq \gamma(\bold{x}) \leq 1$.

To ensure the smoothness of the material layout, we employ convolution filtering~\cite{bourdin2001filters,bruns2001topology} to prevent the occurrence of checkerboard patterns in material distribution~\cite{diaz1995checkerboard}, as follows:
\begin{equation}
\label{filter}
\tilde{\gamma}_{k} = \frac{\sum_{i \in N_{e, k}}w(\bold{x}_{\mathit{i}})\gamma_{\mathit{i}}}{\sum_{i \in N_{e, k}}w(\bold{x}_{\mathit{i}})},
\end{equation}
where $\tilde{\gamma}_{k}$ is the filtered design variable on node $k$, $N_{e, k}$ is the set of nodes within the filtering radius $R$, and $ w(\bold{x}_{\mathit{i}})$ is the weighting function for filtering defined in the following equation.
\begin{equation}
w(\bold{x}_{\mathit{i}}) = \mathit{R} - \lvert \bold{x}_{\mathit{i}} - \bold{x}_{\mathit{k}}\rvert.
\end{equation}

Although the filter can ensure a smooth transition of design variables, it also results in more transitional regions with intermediate values between 0 and 1. Therefore, a projection function is used to eliminate transitional grey values, as follows~\cite{wang2011projection}:
\begin{equation}
\label{projection}
\bar{\tilde{\gamma}}_{k} = \frac{\tanh(\beta\eta) + \tanh(\beta(\tilde{\gamma}_{k} - \eta))}{\tanh(\beta\eta) + \tanh(\beta(1 - \eta))},
\end{equation}
where $\bar{\tilde{\gamma}}_{k}$ is the design variable after projection; $\beta$ is the projection parameter, which controls the sharpness of contours; and $\eta$ is the threshold parameter.

\subsection{Fluid and solid representation using the design variable field}
\label{Fluid and solid representation using the design variable field}
As mentioned in section~\ref{Basic concept of topology optimization}, we employ a relaxed characteristic function for topology optimization on particle flow. In this study, $\gamma = 1$ 
represents the fluid region in $\mathit{\Omega}$, and $\gamma = 0$ represents the pseudo-solid region in $D \backslash \mathit{\Omega}$. To allocate fluid and solid in the design domain, according to the previous research for fluid topology optimization~\cite{borrvall2003topology}, the fictitious body force is required to incorporate into the momentum equations~(\ref{gef}) and~(\ref{gep}) for both phases, which are re-formulated as follows:
\begin{align}
&\nabla \cdot ({\phi_{\rm f} \rho_{\rm f}\bold{u}_{\rm f} \bold{u}_{\rm f}}) =  -\phi_{\rm f} \nabla p + \nabla \cdot (\mu \phi_{\rm f}(\nabla \bold{u}_{\rm f} + \nabla \bold{u}^\intercal_{\rm f})) + \phi_{\rm f} \rho_{\rm f} \bold{g} - \phi_{\rm f}\alpha_{\rm f} \bold{u}_{\rm f}, \\
&\nabla \cdot ({\phi_{\rm p} \rho_{\rm p} \bold{u}_{\rm p} \bold{u}_{\rm p}}) =  -\phi_{\rm p} \nabla p + \phi_{\rm p} \rho_{\rm p} \bold{g} + K(\bold{u_{\rm f}} - \bold{u_{\rm p}})- \phi_{\rm p}\alpha_{\rm p} \bold{u}_{\rm p}. 
\end{align}
The inverse permeability $\alpha$ can be expressed as follows:
\begin{equation}
\alpha_{l}(\bar{\tilde{\gamma}}) = \overline{\alpha}_{l} + (\underline{\alpha_{l}} - \overline{\alpha}_{l})\bar{\tilde{\gamma}} \frac{1 + q}{\bar{\tilde{\gamma}} + q},
\end{equation}
where $l$ denotes the phase with $l \in \{ \rm{f},\rm{p}\}$, $\overline{\alpha}_{l}$ and $\underline{\alpha_{l}}$ are the maximum and minimum values of the inverse permeability, respectively, and $q$ is the convex parameter used to control the degree of penalization for intermediate media. It is crucial to note that since the inertia of particles is much greater than that of the fluid, using the same maximum value of the fictitious body force may lead to undesirable penetration of particles into solid regions. To ensure a more realistic simulation, we adjust the maximum fictitious body force for particles as follows, taking into account their greater inertia:
\begin{equation}
\overline{\alpha}_{\rm p} = c_{\rm \alpha}\overline{\alpha}_{\rm f},
\end{equation}
where $\overline{\alpha}_{\rm p}$ is the maximum value of the inverse permeability for the particle phase, $\overline{\alpha}_{\rm f}$ is the maximum value of the inverse permeability for the fluid phase, and $c_{\alpha}$ is the coefficient to increase the fictitious body force for the particle phase. Based on a sufficiently large value of the fictitious body force, it is difficult for both phases to pass through the solid region ($\gamma = 0$). Note that though the issue of particle penetration can be addressed by adjusting the maximum value of the fictitious body force, the shear stress exerted by the pseudo-solids on the particles is simplified to no-slip conditions due to a limitation in the mathematical representation of Darcy's force. Despite this simplification, sensitivities can be utilized to obtain the optimized distribution of pseudo-solids using the design variable field.

\subsection{Maximization problem: particle drag variation}
Under fixed operating conditions, increasing particle resistance is typically achieved through flow field design that either lengthens the particle path or introduces obstacles within the device. Due to the relatively longer relaxation time, particles cannot instantly follow changes in fluid velocity, particularly when encountering obstacles. This results in a velocity difference between the particles and the fluid. Additionally, the continuous changes in the flow path prevent the particles from reaching their terminal velocity, further sustaining the velocity difference. This sustained velocity difference enhances the drag force acting on the particles, as described by the term $K(\bold{u_{\rm f}} - \bold{u_{\rm p}})$ in Eq.~(\ref{gep}). As the particles undergo repeated acceleration and deceleration due to obstacles or turns in the flow path, the drag force acting on them varies accordingly. Based on this, the variation in particle drag force can be regarded as a key indicator of the resistance experienced by the particles. Therefore, we formulate a topology optimization problem to maximize the magnitude of particle drag variation, while also considering power dissipation to meet industrial application requirements. The expression of the optimization problem is shown below:
\begin{equation}
\begin{split}
&\text{maximize }   J = w_{\zeta} \frac{\zeta}{\zeta_{0}} - w_{\Phi}\frac{\Phi}{\Phi_{0}}, \\
&\begin{split}\text{subject to } V =&\int_{\mathit{\Omega}}\gamma(\bold{x})d\mathit{\Omega} - \overline{V} \leq 0, \\
&0 \leq \gamma(\bold{x}) \leq 1 \text{ for } \forall \bold{x} \in D,\end{split}
\end{split}
\end{equation}
where $J$ is the multi-objective function, $w_{\Phi}$ and $w_{\zeta}$ are the weighting factors, $\Phi$ is the power dissipation of the fluid phase, $\zeta$ is a function related to the variation of drag force on the particle phase, and $\overline{V}$ is the maximum value of the volume constraint. Due to the difference in scale between $\Phi$ and $\zeta$, they are normalized based on their values $\Phi_{0}$ and $\zeta_{0}$ from the first iteration of the optimization process. The functions $\Phi$ and $\zeta$ can be expressed as follows:
\begin{align}
&\Phi = \int_{\mathit{\Omega}} \left[ \frac{1}{2}\mu(\nabla \bold{u}_{\rm f} + \nabla \bold{u}^\intercal_{\rm f}) : (\nabla \bold{u}_{\rm f} + \nabla \bold{u}^\intercal_{\rm f})+\alpha(\bar{\tilde{\gamma}}) {\lvert\bold{u}_{\rm f}\rvert}^2 \right] d\mathit{\Omega}, \\[10pt]
&\zeta =  \int_{\mathit{\Omega}} (\zeta_{\rm D} - \overline{\zeta}_{\rm D}) d\mathit{\Omega},
\end{align}
where $\zeta_{\rm D}$ represents the magnitude of particle drag force, given by the equation $\zeta_{\rm D} = K \lvert \bold{u}_{\rm f} - \bold{u}_{\rm p} \rvert$, and $\overline{\zeta}_{\rm D}$ is the average value of the magnitude of particle drag force in $\mathit{\Omega}$. 

\section{Numerical implementation}
\label{Numerical implementation}
\subsection{NAPPLE algorithm}
The finite difference method called the NAPPLE algorithm~\cite{lee1992artificial} is employed for the numerical analysis with in-house C++ code. A distinguishing feature of this algorithm is its ability to define all state variables using a collocated grid system. This approach eliminates the need for a staggered grid to address the checkerboard pressure issue~\cite{prakash1985control,date1996complete}, which makes programming more straightforward and computation more efficient. In this section, we will introduce the concepts of the NAPPLE algorithm and numerical discretization.

\begin{figure}[h]
\centering
\includegraphics[width=0.4\textwidth]{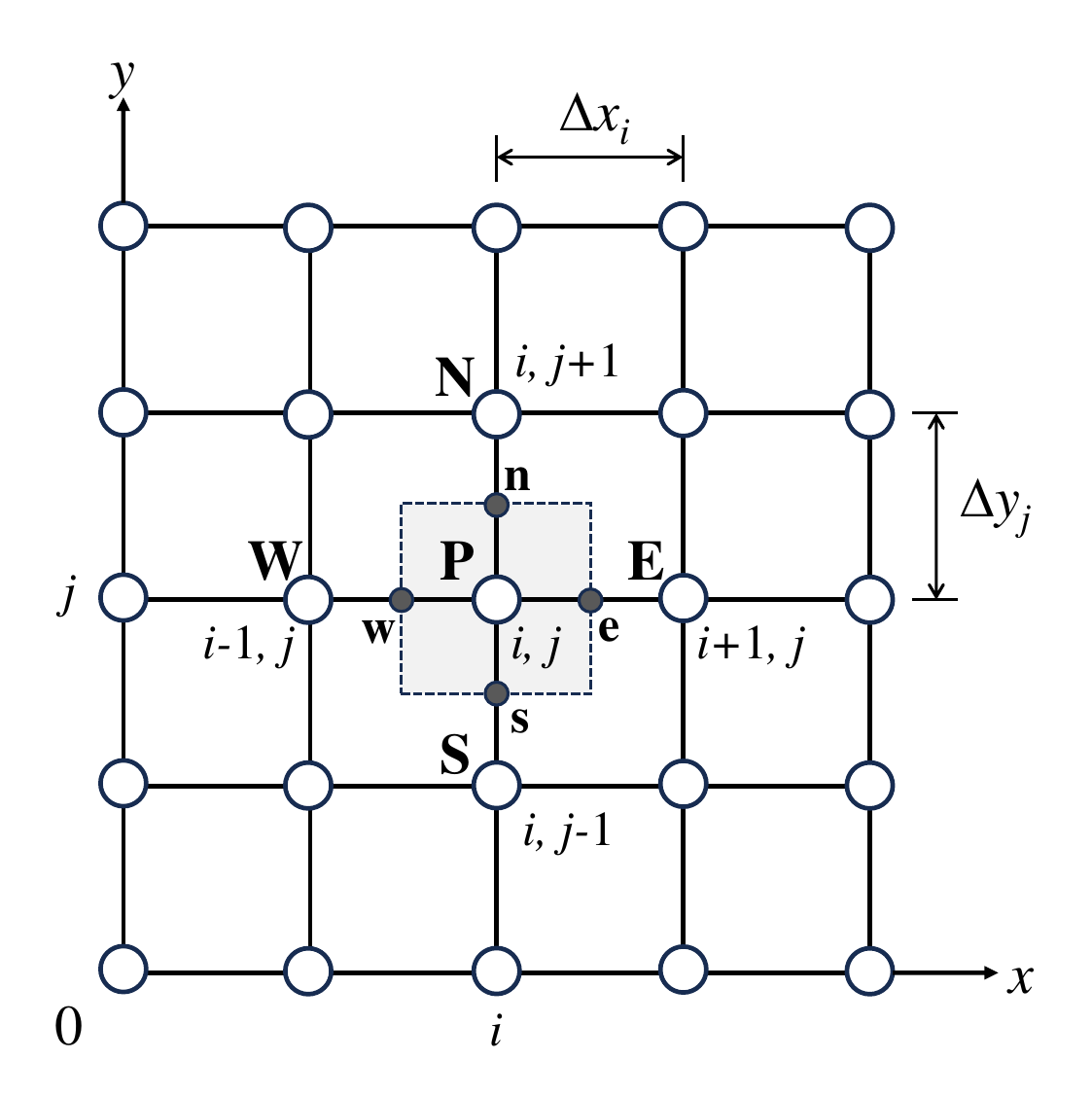} 
\caption{Node Notation for point $\bold{P}$ on a Cartesian grid system.}
\label{Fig.1}
\end{figure}

We consider a uniform Cartesian grid system for spatial discretization. A point $\bold{P}$ is denoted as ${\rm P}(i, j)$, where $i$ and $j$ indicate the column number and the row number of point $\bold{P}$ along the $x$-axis and the $y$-axis, as shown in Fig.~\ref{Fig.1}. To ensure numerical stability, we use the pseudo-transient method to achieve steady-state solutions~\cite{ferziger2019computational}. This method introduces a pseudo-time derivative, allowing the changes in state variables to progressively diminish to zero as pseudo-time steps march forward. Therefore, we consider dimensionless transient governing equations for both phases and rewrite Eqs.~(\ref{gef})--(\ref{pc}) under the assumption of incompressible flow as follows:
\begin{align}
\label{NS11}
&Re\phi_{\rm f}^{*}\frac{\partial \bold{u}_{\rm f}^{*}}{\partial t^{*}} + Re\phi_{\rm f}^{*}(\bold{u}_{\rm f}^{*}\cdot \nabla^{*})\bold{u}_{\rm f}^{*} =  -\phi_{\rm f}^{*} \nabla^{*} \hat{p} + \nabla^{*} \cdot (\phi_{\rm f}^{*}\nabla^{*} \bold{u}_{\rm f}^{*}) + \bold{s}_{\rm f}^{*}, \\[10pt]
\label{ndc1}
&\frac{\partial \phi_{\rm f}^{*}}{\partial t^{*}} + \nabla^{*} \cdot (\phi_{\rm f}^{*} \bold{u}_{\rm f}^{*}) = 0, \\[10pt]
&Re \hat{\rho}\phi_{\rm p}^{*}\frac{\partial \bold{u}_{\rm p}^{*}}{\partial t^{*}} + Re \hat{\rho}\phi_{\rm p}^{*}(\bold{u}_{\rm p}^{*}\cdot \nabla^{*})\bold{u}_{\rm p}^{*} =  -\phi_{\rm p}^{*} \nabla^{*} \hat{p} + \bold{s}_{\rm p}^{*}, \\[10pt]
\label{ndc2}
&\frac{\partial \phi_{\rm p}^{*}}{\partial t^{*}} + \nabla^{*} \cdot (\phi_{\rm p}^{*} \bold{u}_{\rm p}^{*}) = 0,
\end{align}
where $\bold{s}_{\rm f}$ is the source term for the fluid phase, given by $\nabla \cdot (\mu \phi_{\rm f}\nabla \bold{u}^\intercal_{\rm f}) + \phi_{\rm f} \rho_{\rm f} \bold{g}$, and $\bold{s}_{\rm p}$ is the source term for the particle phase, given by $\phi_{\rm p} \rho_{\rm p} \bold{g} + K(\bold{u_{\rm f}} - \bold{u_{\rm p}})$. Notably, since the transpose term in the viscous term will be treated explicitly in the discretization procedure, it is included in the source term. The variables with asterisks represent the dimensionless variables defined as follows:
\begin{equation}
\begin{split}
&\bold{u}_{\rm f}^{*} = \frac{\bold{u}_{\rm f}}{U} \ , \ \bold{u}_{\rm p}^{*} = \frac{\bold{u}_{\rm p}}{U} \ , \ \phi_{\rm f}^{*} = \frac{\phi_{\rm f}}{\phi} \ , \ \phi_{\rm p}^{*} = \frac{\phi_{\rm p}}{\phi} \ , \ \bold{s}_{\rm f}^{*} = \frac{L^{2}}{\mu U}\bold{s}_{\rm f} \  , \ \bold{s}_{\rm p}^{*} = \frac{L^{2}}{\mu U}\bold{s}_{\rm p} \ , \ \\[10pt]
&\nabla^{*} = L\nabla \  ,\ \ t^{*} = \frac{U}{L}t \ , \ Re = \frac{\rho U L}{\mu} \  , \ \ p^{*} = \frac{p - p_{0}}{\rho_{\rm f} U^{2}} \  , \ \hat{p} = Rep^{*}, \ \hat{\rho} = \frac{\rho_{\rm p}}{\rho_{\rm f}}, 
\end{split}
\end{equation}
where $U$ is the characteristic velocity, $L$ is the characteristic length, and $Re$ is the flow Reynolds number. For simplicity, the dimensionless variables used in this section will be written without the asterisk in the following description.

Due to the implicit discretization scheme used in the NAPPLE algorithm, the governing equations are discretized into a matrix of linear algebraic equations to solve for the state variables. For the discretization process, the weighting function scheme~\cite{shong1989weighting} is employed for the momentum equations of the fluid phase, and a well-known upwind scheme~\cite{anderson1995computational} is utilized for the particle phase to ensure numerical stability. Owing to its general familiarity, the explanation of the upwind scheme employed for the particle phase will not be expanded upon in this section. Instead, we will focus on detailing the weighting function scheme used for the fluid phase. Considering a point ${\rm P}(i, j)$ on a structured grid shown in Fig.~\ref{Fig.1}, the discretized momentum equations~(\ref{NS11}) in the $x$-direciton and the $y$-direction can be expressed as follows:
\begin{align}
\label{DNS11}
&(a_{\rm W})_{i,j} u_{i - 1,j}+ (a_{\rm S})_{i,j} u_{i,j - 1} + (a_{\rm P})_{i,j} u_{i,j} 
+ (a_{\rm E})_{i,j} u_{i + 1,j}+ (a_{\rm N})_{i,j} u_{i,j + 1} = (r_{u})_{i,j} + \frac{\partial \hat{p}}{\partial x}, \\[10pt]
\label{DNS2}
&(a_{\rm W})_{i,j} v_{i - 1,j}+ (a_{\rm S})_{i,j} v_{i,j - 1} + (a_{\rm P})_{i,j} v_{i,j} 
+ (a_{\rm E})_{i,j} v_{i + 1,j}+ (a_{\rm N})_{i,j} v_{i,j + 1} = (r_{v})_{i,j} + \frac{\partial \hat{p}}{\partial y},
\end{align}
where $u$ is the $x$-component of velocity, $v$ is the $y$-component of velocity, $a$ is the weighting coefficient and the subscripts $\rm W, E, S, N$ denote the neighborhood grid points of point $\bold{P}$ as shown in Fig.~\ref{Fig.1}. The terms $r_{u}$ and $r_{v}$ are the remainders in Eqs.~(\ref{DNS11}) and~(\ref{DNS2}). The weighting coefficients can be further expressed with the weighting function scheme as follows:
\begin{equation}
\label{wfs}
\begin{split}
&(a_{\rm W})_{i,j} = \frac{w_{f}((\mathit{z}_{\rm w})_{\mathit{i,j}})}{\Delta x^{2}}, 
  (a_{\rm E})_{i,j} = \frac{w_{f}(-(\mathit{z}_{\rm e})_{\mathit{i,j}})}{\Delta x^{2}}, 
  (a_{\rm S})_{i,j} = \frac{w_{f}((\mathit{z}_{\rm s})_{\mathit{i,j}})}{\Delta y^{2}}, 
  (a_{\rm N})_{i,j} = \frac{w_{f}(-(\mathit{z}_{\rm n})_{\mathit{i,j}})}{\Delta y^{2}}, \\[10pt]
&(a_{\rm P})_{i,j} = -(a_{\rm W})_{i,j} - (a_{\rm S})_{i,j} - (a_{\rm E})_{i,j} - (a_{\rm N})_{i,j} - \frac{Re}{\Delta t}, 
\quad (r_{u})_{i,j} = -\left( \frac{s_{x}}{\phi} +\frac{Re\ u_{0}}{\Delta t}\right)_{\rm f},
\quad (r_{v})_{i,j} = -\left( \frac{s_{y}}{\phi} +\frac{Re\ v_{0}}{\Delta t}\right)_{\rm f},
\end{split}
\end{equation}
where the subscript $\rm f$ denotes the fluid phase, the subscripts $\rm w, e, s, n$ denote the midpoint between each adjacent interval of point $\bold{P}$, $u_{0}$ is the $x$-component of velocity at the previous time step, and $v_{0}$ is the $y$-component of velocity at the previous time step. The weighting function $w_{f}(\mathit{z})$ and the parameter $z$ used in Eq.~(\ref{wfs}) are defined as follows:
\begin{equation}
\begin{split}
\label{wf}
&w_{f}(\mathit{z})  =\frac{\mathit{z}}{1-e^{-\mathit{z}}},			\\[5pt]
&(z_{\rm w})_{i,j} = \left( Re \overline{u}_{i-1} - \frac{1}{\overline{\phi}_{i-1}}\frac{\partial \phi}{\partial x}\right)_{\rm f} \Delta x, \quad 
(\mathit{z}_{\rm e})_{i,j} = \left( Re \overline{u}_{i} - \frac{1}{\overline{\phi}_{i}}\frac{\partial \phi}{\partial x}\right)_{\rm f} \Delta x, \\[5pt]
&(z_{\rm s})_{i,j} =\left( Re \overline{v}_{j-1} - \frac{1}{\overline{\phi}_{j-1}}\frac{\partial \phi}{\partial y}\right)_{\rm f} \Delta y, \quad
(\mathit{z}_{\rm n})_{i,j} = \left( Re \overline{v}_{j} - \frac{1}{\overline{\phi}_{j}}\frac{\partial \phi}{\partial y}\right)_{\rm f} \Delta y,
\end{split}
\end{equation}
where $\overline{u}_{i}$, $\overline{v}_{j}$,  $\overline{\phi}_{i}$ and $\overline{\phi}_{j}$ represent the average velocity and the average volume fraction in the intervals $x_{i-1} \leq x \leq x_{i}$ and $y_{i-1} \leq y \leq y_{i}$, respectively. For the weighting function $w_{f}(\mathit{z})$ in Eq.~(\ref{wf}), we use the power-law approximation to reduce the computational cost~\cite{patankar2018numerical}, given by the equation $w_{f}(\mathit{z})  \approx [ 0, (1 - 0.1\lvert z \rvert)^{5} ] +  [0, z]$. Besides, as suggested by Lee~\cite{lee1988new}, the unsteady term in Eq.~(\ref{NS11}) is discretized with implicit form as follows:
\begin{equation}
\frac{\partial \bold{u}_{\rm f}}{\partial t} = \left( \frac{\bold{u} - \bold{u}_{\tiny 0}}{\Delta t} \right)_{\rm f}.
\end{equation}
Though Eqs.~(\ref{DNS11}) and~(\ref{DNS2}) are solvable with the initial guess of the pressure field, it is not guaranteed to satisfy the continuity equation. Therefore, it is required to update state variables iteratively, similar to the well-known SIMPLE algorithm~\cite{patankar1983calculation}. Unlike the SIMPLE algorithm, which relies on a pressure-correction equation, the NAPPLE algorithm integrates the discretized continuity equation and momentum equations to derive a pressure-linked equation. This equation is then directly used to compute the pressure field. To derive the pressure-linked equation, we first recast Eqs.~(\ref{DNS11}) and~(\ref{DNS2}) into a compact form for the volume flux:
\begin{align}
\label{CDNS11}
&(\phi_{i, j}u_{i, j})_{\rm f} = (\phi_{i, j}\hat{u}_{i, j})_{\rm f} - \kappa_{i,j}\left(\frac{ \partial \hat{p}}{\partial x}\right),	\\[10pt]
\label{CDNS111}
&(\phi_{i, j}v_{i, j})_{\rm f} = (\phi_{i, j}\hat{v}_{i, j})_{\rm f}  - \kappa_{i,j}\left(\frac{ \partial \hat{p}}{\partial y}\right),	\\[10pt]
\label{CDNS2}
&\hat{u}_{i, j} = ((r_{u})_{i,j} - (a_{\rm W})_{i,j} u_{i - 1,j} - (a_{\rm S})_{i,j} u_{i,j - 1} - (a_{\rm P})_{i,j} u_{i,j} - (a_{\rm E})_{i,j} u_{i + 1,j} -  (a_{\rm N})_{i,j} u_{i,j + 1}) / (a_{\rm P})_{i,j},\\[10pt]
\label{CDNS22}
&\hat{v}_{i, j} = ((r_{v})_{i,j} - (a_{\rm W})_{i,j} v_{i - 1,j} - (a_{\rm S})_{i,j} v_{i,j - 1} - (a_{\rm P})_{i,j} v_{i,j} - (a_{\rm E})_{i,j} v_{i + 1,j} -  (a_{\rm N})_{i,j} v_{i,j + 1}) / (a_{\rm P})_{i,j},\\[10pt]
&\kappa_{i,j} = - \left( \frac{{\phi}}{ (a_{\rm P})_{i,j}} \right)_{\rm f}.
\end{align}
As suggested by Moukalled et al.~\cite{moukalled2003pressure}, we sum the continuity equations~(\ref{ndc1}) and~(\ref{ndc2}) to obtain the global continuity equation. Subsequently, the global continuity equation can be utilized to derive the pressure-linked equation, thereby solving for the pressure term shared by both phases. The discretized global continuity equation for point ${\rm P}(i, j)$ can be expressed in the following form:
\begin{equation}
\label{globalcontinuity}
\sum_{l \in \{ \rm f, p \} } \left( \frac{\phi_{i+1, j}u_{i+1, j} - \phi_{i-1, j}u_{i-1, j}}{2 \Delta x} + \frac{\phi_{i, j+1}v_{i, j+1} - \phi_{i, j-1}v_{i, j-1}}{2 \Delta y} \right)_{l} = 0.
\end{equation}
It should be noted that since the sum of the volume fraction of both phases is equal to one at any given time, the influence of the transient term in the global continuity equation can be neglected. Furthermore, as suggested by the prior work~\cite{lee1992artificial}, we adopted the assumptions shown below during the discretization process: 
\begin{align}
\label{assumption1}
&\frac{\kappa_{i+1, j} \left(\frac{\partial \hat{p}}{\partial x}\right)_{i+1, j} - \kappa_{i-1, j} \left(\frac{\partial \hat{p}}{\partial x}\right)_{i-1, j}}{2 \Delta x} \approx \frac{\partial}{\partial x}\left(\kappa \frac{\partial \hat{p}}{\partial x}\right),	\\[10pt]
\label{assumption2}
&\frac{\kappa_{i, j+1} \left(\frac{\partial \hat{p}}{\partial y}\right)_{i, j+1} - \kappa_{i, j-1} \left(\frac{\partial \hat{p}}{\partial y}\right)_{i, j-1}}{2 \Delta y} \approx \frac{\partial}{\partial y}\left(\kappa \frac{\partial \hat{p}}{\partial y}\right).
\end{align}
Based on the assumptions in Eqs.~(\ref{assumption1}) and~(\ref{assumption2}),  integrating Eqs.~(\ref{CDNS11})--(\ref{globalcontinuity}) yields the pressure-linked equation, which eliminates the checkerboard pressure pattern~\cite{lee1992artificial}. The pressure-linked equation is formulated as follows:
\begin{align}
\label{pleq}
&\frac{\partial}{\partial x}\left(\kappa \frac{\partial \hat{p}}{\partial x}\right) + \frac{\partial}{\partial y}\left(\kappa \frac{\partial \hat{p}}{\partial y}\right) = \hat{\varepsilon}, \\[10pt]
\label{dilationCalc}
&\hat{\varepsilon}_{i, j} = \sum_{l \in \{ \rm f, p \}} \left(\frac{\phi_{i+1, j}\hat{u}_{i+1, j} - \phi_{i-1, j}\hat{u}_{i-1, j}}{2 \Delta x} + \frac{\phi_{i, j+1}\hat{v}_{i, j+1} - \phi_{i, j-1}\hat{v}_{i, j-1}}{2 \Delta y}\right)_{l},
\end{align}
where $\hat{\varepsilon}$ is the pseudo-dilation. The pressure-linked equation can also be discretized into a linear algebraic equation with the harmonic scheme~\cite{patankar2018numerical}:
\begin{equation}
\label{peq}
(b_{\rm W})_{i,j} p_{i - 1,j}+ (b_{\rm S})_{i,j} p_{i,j - 1} + (b_{\rm P})_{i,j} p_{i,j} + (b_{\rm E})_{i,j} p_{i + 1,j}+ (b_{\rm N})_{i,j} p_{i,j + 1} = (r_{p})_{i,j}.
\end{equation}
The coefficient $b$ in Eq.~(\ref{peq}) can be further expressed as follows:
\begin{equation}
\begin{split}
&(b_{\rm W})_{i,j} = \frac{(\kappa_{\rm w})_{i,j}}{\Delta x^{2}}, 
(b_{\rm E})_{i,j} = \frac{(\kappa_{\rm e})_{i,j}}{\Delta x^{2}}, 
(b_{\rm S})_{i,j} = \frac{(\kappa_{\rm s})_{i,j}}{\Delta y^{2}}, 
(b_{\rm N})_{i,j} = \frac{(\kappa_{\rm n})_{i,j}}{\Delta y^{2}}, \\[10pt]
&(b_{\rm P})_{i,j} = -(b_{\rm W})_{i,j} - (b_{\rm S})_{i,j} - (b_{\rm E})_{i,j} - (b_{\rm N})_{i,j}, \ (r_{p})_{i,j} = \hat{\varepsilon}_{i,j},
\end{split}
\end{equation}
Since the global pressure is shared by both the fluid phase and the particle phase, the total pressure conductivity $\kappa$ represents the sum of the pressure conductivities of each phase. In each phase, the pressure conductivity can be obtained through the coefficients in their respective momentum equations. Therefore, the total pressure conductivity $\kappa$ can be expressed as follows:
\begin{align}
&(\kappa_{\rm w})_{i,j}=\sum_{l \in \{ \rm f, p \}}\left(\frac{2}{(\kappa_{i,j})^{-1} + (\kappa_{i-1,j})^{-1}}\right)_{l}
			       =\sum_{l \in \{ \rm f, p \}}\left(\frac{-2\phi_{i,j}\phi_{i-1,j}}{(a_{\rm P})_{i,j}\phi_{i-1,j} + \phi_{i,j}(a_{\rm P})_{i-1,j}}\right)_{l}	,\\[10pt]
&(\kappa_{\rm e})_{i,j}=\sum_{l \in \{ \rm f, p \}}\left(\frac{2}{(\kappa_{i,j})^{-1} + (\kappa_{i+1,j})^{-1}}\right)_{l}
			       =\sum_{l \in \{ \rm f, p \}}\left(\frac{-2\phi_{i,j}\phi_{i+1,j}}{(a_{\rm P})_{i,j}\phi_{i+1,j} + \phi_{i,j}(a_{\rm P})_{i+1,j}}\right)_{l},	\\[10pt]
&(\kappa_{\rm s})_{i,j}=\sum_{l \in \{ \rm f, p \}}\left(\frac{2}{(\kappa_{i,j})^{-1} + (\kappa_{i,j-1})^{-1}}\right)_{l}
                                  =\sum_{l \in \{ \rm f, p \}}\left(\frac{-2\phi_{i,j}\phi_{i,j-1}}{(a_{\rm P})_{i,j}\phi_{i,j-1} + \phi_{i,j}(a_{\rm P})_{i,j-1}}\right)_{l},	\\[10pt]
&(\kappa_{\rm n})_{i,j}=\sum_{l \in \{ \rm f, p \}}\left(\frac{2}{(\kappa_{i,j})^{-1} + (\kappa_{i,j+1})^{-1}}\right)_{l}
			       =\sum_{l \in \{ \rm f, p \}}\left(\frac{-2\phi_{i,j}\phi_{i,j+1}}{(a_{\rm P})_{i,j}\phi_{i,j+1} + \phi_{i,j}(a_{\rm P})_{i,j+1}}\right)_{l}.
\end{align}
\begin{algorithm}[h]
\caption{NAPPLE algorithm for solving particle flow problem}\label{alg:napple}
\begin{algorithmic}[1]
\State All the state variables and the parameters are initialized.
\While{the solution of state variables does not converge}
\For{each phase $l$ $\in \{ \rm f, p\}$}
\State Impose boundary conditions in Eqs.~(\ref{bc1}) and (\ref{bc12}) for velocities $u_{l}, v_{l}$
\EndFor
\State Impose the boundary condition in Eq.~(\ref{bcp}) for the global pressure $\hat{p}$ shared by both phases
\State Impose the boundary condition in Eqs.~(\ref{bcp2}) and (\ref{bc2}) for the particle volume fraction $\phi_{\rm p}$
\For{each phase $l$ $\in \{ \rm f, p\}$}
\State Compute the weighting factors $(a_{\rm W})_{l}, (a_{\rm E})_{l}, (a_{\rm S})_{l}, (a_{\rm N})_{l}, (a_{\rm P})_{l}$ in Eqs.~(\ref{DNS11}) and~(\ref{DNS2})
\State Compute the pseudo-velocities $\hat{u}_{l}$, $\hat{v}_{l}$ in Eqs.~(\ref{CDNS2}) and~(\ref{CDNS22})
\EndFor
\State Compute the weighting factors $b_{\rm W}, b_{\rm E}, b_{\rm S}, b_{\rm N}, b_{\rm P}$ in the left-hand side of Eq.~(\ref{peq})
\State Compute the pseudo-dilation $\hat{\varepsilon}$ in the right-hand side of Eq.~(\ref{peq}) for both phases, which expressed by Eq.~(\ref{dilationCalc})
\State Solve the matrix of the pressure-linked equation~(\ref{peq}) with a matrix solver
\State Update the pressure $\hat{p}$ with a relaxation factor $\omega_{\rm p}$, $\hat{p} \gets (1 - \omega_{\rm p})\hat{p}^{k-1} + \omega_{\rm p}\hat{p}^{k}$
\State Compute the pressure gradients $\frac{\partial \hat{p}}{\partial x}$ and $\frac{\partial \hat{p}}{\partial y}$
\For{each phase $l$ $\in \{ \rm f, p\}$}
\State Add the pressure gradients $\frac{\partial \hat{p}}{\partial x}$ and $\frac{\partial \hat{p}}{\partial y}$ to the right-hand side of momentum equations
\State Solve the matrix of momentum equations~(\ref{DNS11}) and~(\ref{DNS2})with a matrix solver
\State Update the velocity $u_l$ with a relaxation factor $\omega_{\rm v}$, $u_{l} \gets (1 - \omega_{\rm v})u^{k-1}_{l} + \omega_{\rm v} u^{k}_{l}$
\State Update the velocity $v_l$ with a relaxation factor $\omega_{\rm v}$, $v_{l} \gets (1 - \omega_{\rm v})v^{k-1}_{l} + \omega_{\rm v} v^{k}_{l}$
\EndFor
\State Discretize the continuity equation for the particle phase in Eq.~(\ref{ndc2}) with a numerical scheme
\State Obtain the particle volume fraction $\phi_{\rm p}$ by solving the matrix of the continuity equation~(\ref{ndc2})
\State Obtain $\phi_{\rm f}$ with the relationship $\phi_{\rm f} + \phi_{\rm p} = 1$
\State Update the iteration number $k$, $k \gets k+1$
\EndWhile
\State Return the final solution of state variables of each phase
\end{algorithmic}
\end{algorithm}

Based on the description provided, the steps of the NAPPLE algorithm for particle flow can be summarized as pseudocode in Algorithm~\ref{alg:napple}. For particle flow analysis, we must determine the particle volume fraction by discretizing the particle continuity equation~(\ref{ndc2}) to obtain the linear algebraic equation, which is then solved to find the volume fraction values. The upwind scheme is utilized to discretize the continuity equation of the particle phase.

Besides, the SIS algorithm (Strongly implicit solver) suggested by Lee~\cite{lee1990strongly} is utilized for solving the matrices of linear algebraic equations. To ensure numerical stability, the relaxation factors $\omega_{\rm p}$ and $\omega_{\rm v}$ are utilized to update the pressure and velocities for both phases and the values of $\omega_{\rm p}$ and $\omega_{\rm v}$ are set to 0.2 and 0.4 in this study. Moreover, we must iterate the numerical procedure until we achieve solution convergence. The convergence criteria are as follows:
\begin{align}
&\mathop{\max}_{{\rm P}} \left| u^{k}_{i, j} - u^{k - 1}_{i, j} \right| \leq \rm{tol}_\mathit{u}, \\
&\mathop{\max}_{{\rm P}} \left| v^{k}_{i, j} - v^{k - 1}_{i, j} \right| \leq \rm{tol}_\mathit{v},
\end{align}
where $\rm{tol}_\mathit{u}$ and $\rm{tol}_\mathit{v}$ represent the convergence criteria for the velocity component $u$ and $v$ for both phases and are set to $1 \times 10^{-5}$, respectively.

\subsection{Optimization algorithm}
\begin{figure}[h]
\centering
\includegraphics[width=1.0\textwidth]{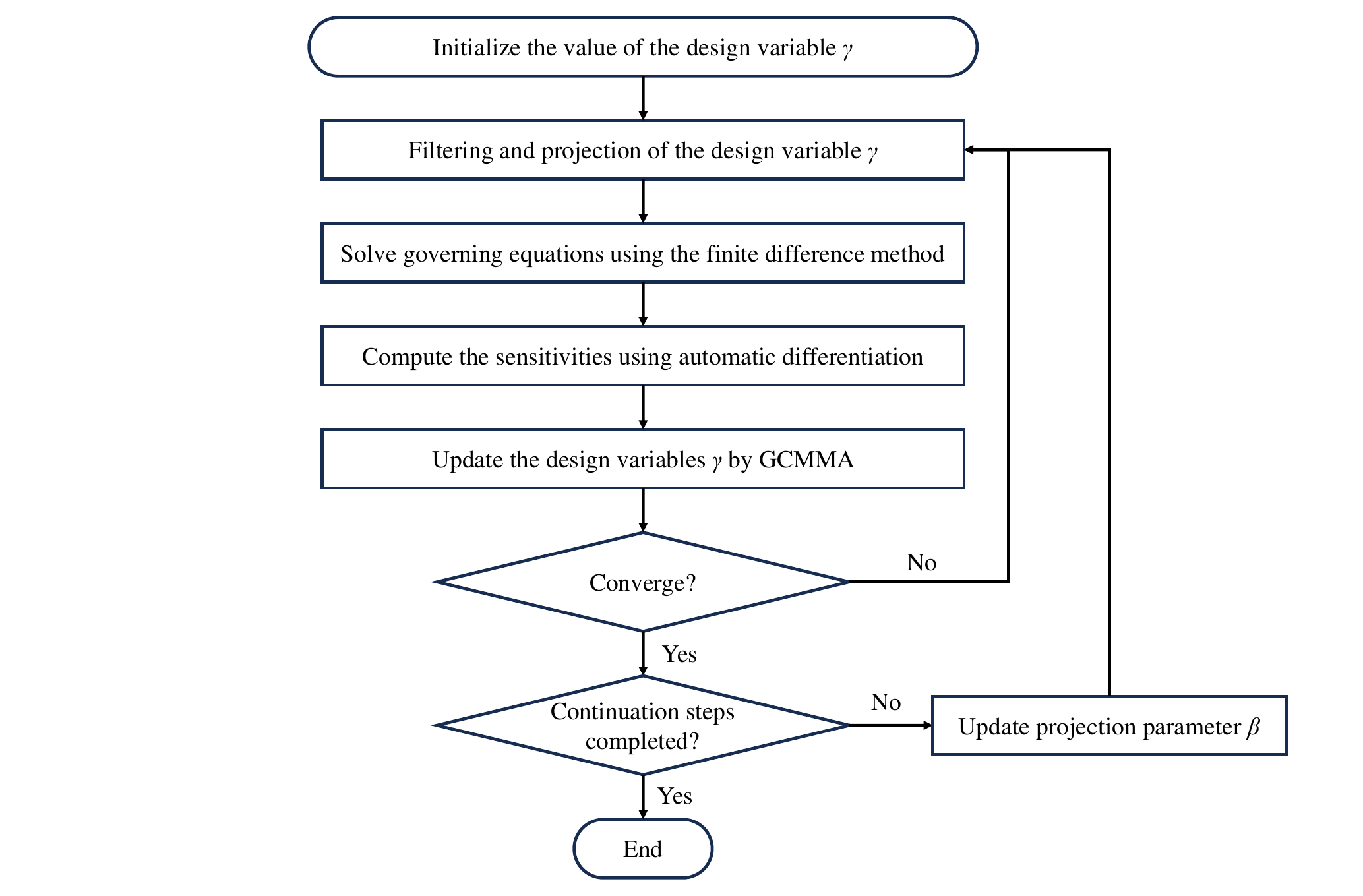} 
\caption{Flowchart of topology optimization process with a continuation approach of the projection parameter $\beta$.}
\label{Fig.2}
\end{figure}

In the optimization process, we use a continuation method to progressively increase the projection parameter $\beta$ in Eq.~(\ref{projection}) in each continuation step, which helps to eliminate intermediate design variables. The optimization algorithm iteratively solves the problem based on gradient information during each of these steps. The flowchart of the topology optimization process is shown in Fig.~\ref{Fig.2}.

As shown in Fig.~\ref{Fig.2}, the design variables $\gamma$ are initialized in the first step. In the second step, the design variables are regularized with the convolution filter in Eq.~(\ref{filter}) and the projection function in Eq.~(\ref{projection}). We solve the governing equations using the finite difference method in the third step and compute the sensitivities in the fourth step. In the fifth step, we use the GCMMA algorithm to update the design variables~\cite{CCSA}. If the objective function does not meet the convergence criterion, the procedure returns to the second step. If it meets the criterion, the procedure verifies whether the continuation steps are completed. If not, the projection parameter $\beta$ in Eq.~(\ref{projection}) is doubled in the next continuation step, and the procedure returns to the second step. Otherwise, the optimization procedure is terminated. The convergence criterion for the objective function is defined as follows:
\begin{equation}
\left| \frac{ J_{n} - J_{n-1}}{J_{n}} \right| < \epsilon,
\end{equation}
where $\epsilon$ is the criterion for the optimization convergence and is set to $1 \times 10^{-4}$ in this study.
\subsection{Sensitivity analysis based on automatic differentiation}
\label{Sensitivity analysis based on automatic differentiation}
Automatic differentiation offers a robust method for calculating derivatives of functions without the analytical adjoint formulations. Composite functions are decomposed into elementary operations in the process of automatic differentiation. By applying the chain rule to each operation, this approach enables precise computation of derivatives. Therefore, automatic differentiation is well-suited for handling the relatively complex equations for particle flow in this study. As reverse differentiation traces the computation backward from the outputs to the inputs, it solves gradient-based optimization problems involving functions with multiple inputs effectively. In this study, we implement the reverse mode automatic differentiation using the library Adept C++~\cite{hogan2014fast}. For sensitivity analysis, it is necessary to declare state variables and related variables as \texttt{adouble} objects defined in Adept C++, enabling the tracking of all calculations and the reverse computation of derivatives. To verify the correctness of the sensitivity calculations, the sensitivities obtained from Apept C++ are compared with the results from a finite difference check, as shown in~\ref{sensitivity appendix}.

Similar to the iterative solution-based study by Towara and Naumann~\cite{towara2013discrete}, since the reverse mode automatic differentiation requires recording computational details from each iteration for sensitivity analysis, it significantly increases the memory requirement for the iterative process. To address this issue, we employ a checkpointing algorithm for the sensitivity analysis~\cite{griewank2008evaluating,griewank2000algorithm}. The operation of a checkpointing algorithm involves dividing the iterative process into several smaller subintervals. At each checkpoint, state variables are recorded and subsequently used to compute adjoint information in reverse, which serves as the initial condition for the next subinterval. As a result, this method significantly reduces memory requirements due to the storage of essential data at these checkpoints and the adjoint calculations within each subinterval.
\section{Results and discussion}
\label{Results and discussion}
\subsection{Validation of numerical scheme }
\begin{figure*}[h]
\centering
\includegraphics[width=1.0\textwidth]{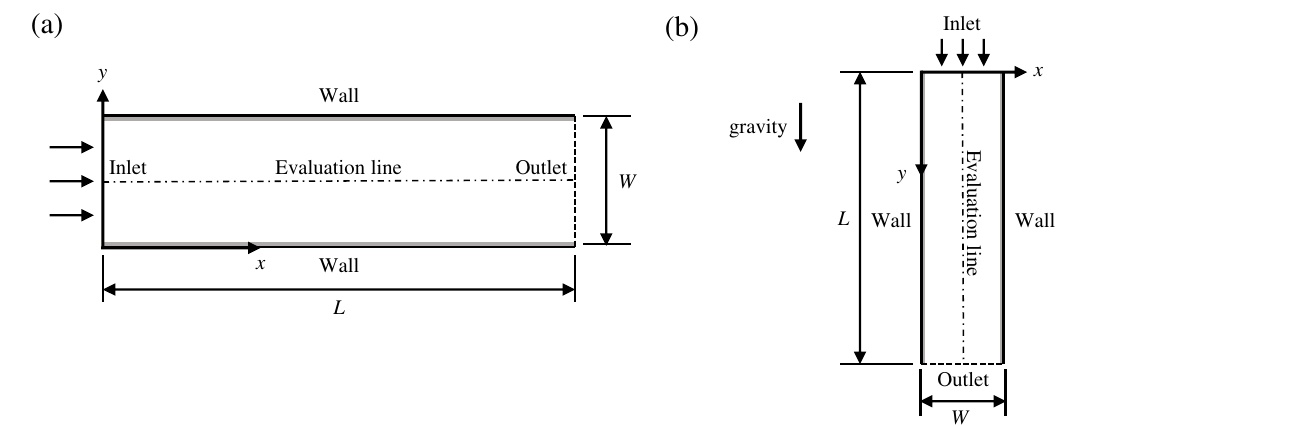} 
\caption{
Dimensions and boundary expressions of analysis domain for numerical validation: (a) Particles flow through a horizontal pipe in the numerical example 1 (b) Particles driven by the gravity effect flow through a vertical pipe in the numerical example 2.}
\label{Fig.15}
\end{figure*}

\begin{figure*}[h]
\centering
\includegraphics[width=0.9\textwidth]{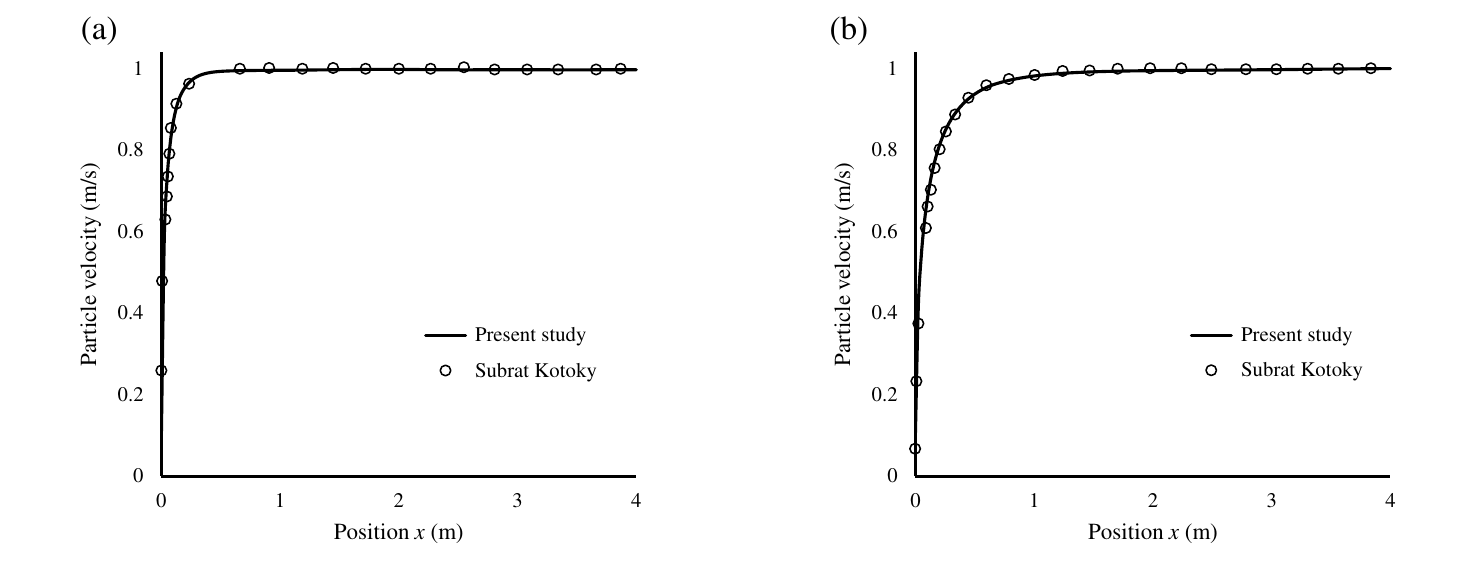} 
\caption{
Comparison of particle velocity distribution along the central line in a horizontal pipe with $\rm \rho_p$ = 665 $\rm kg/m^3$: (a) $ d_{\rm p}$ = 200 $\rm{\mu} m$ (b) $d_{\rm p}$ = 400 $\rm{\mu} m$.}
\label{Fig.3}
\end{figure*}

Using two numerical examples from prior studies~\cite{kotoky2018parametric,kotoky2018development}, we validate the feasibility of our numerical scheme. The fluid domain is discretized with a quadrilateral mesh of $160 \times 20$. In the first numerical example, the velocity distribution on the central line in a horizontal straight pipe serves as the crucial point, as shown in Fig.~\ref{Fig.15}(a). The fluid phase and particle phase pass through the pipe with a velocity of 1.0 m/s and 0.1 m/s, respectively. Free-slip boundary conditions are applied on the walls. Table~\ref{validation1} shows the relevant physical parameters for this setup. We examine the velocity distribution along the central line of the pipe. As shown in Fig.~\ref{Fig.3}, the numerical results closely match the prior studies. With the increases in particle diameter, a longer distance is required to reach the terminal velocity of particles reasonably.

\begin{figure*}[h]
\centering
\includegraphics[width=0.9\textwidth]{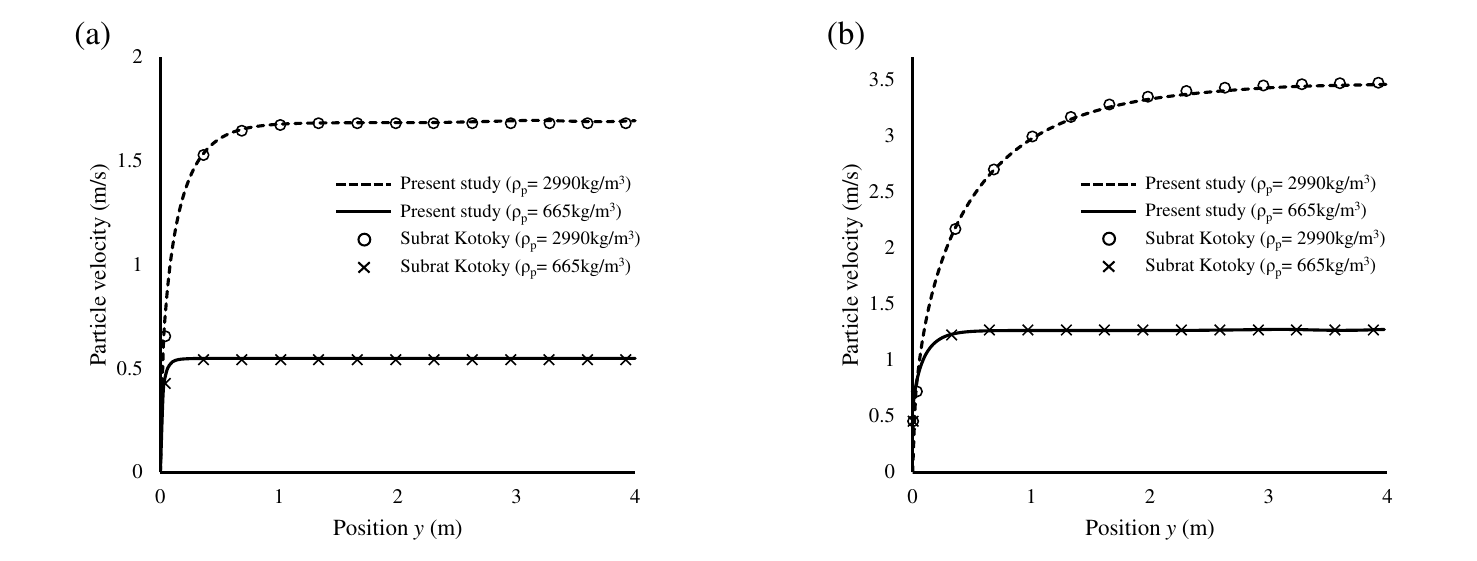} 
\caption{
Comparison of particle velocity distribution along the central line in a vertical pipe: (a) $ d_{\rm p}$ = 200 $\rm{\mu} m$ (b) $d_{\rm p}$ = 400 $\rm{\mu} m$.}
\label{Fig.4}
\end{figure*}

To further validate the flow solver, we consider the effect of gravity. In the second numerical example, the particles driven by the gravity effect flow through a vertical pipe in the quiescent medium, as shown in Fig.~\ref{Fig.15}(b). The velocities of both phases are 0.0001 m/s at the inlet. Free-slip boundary conditions are also applied on the walls. Based on the physical parameters shown in Table~\ref{validation1},  the velocity distribution along the central line remains consistent with the prior studies in Fig.~\ref{Fig.4}. 

\begin{table}[H]\centering
\caption{Physical parameters and dimensions for numerical validation.}
\label{validation1}
\fontsize{1}{1}\selectfont
\resizebox{1.0\linewidth}{!}{
\begin{tabular}{@{}lccll@{}}
																							 \specialrule{.1em}{.05em}{.05em} 
Parameters			& \multicolumn{1}{l}{Symbols} 	& \multicolumn{2}{c}{Values}   		& Unit 			\rule[1.2ex]{0pt}{1.2ex}\\ \specialrule{.1em}{.05em}{.05em} 
Pipe length                   	& $L$                                      & \multicolumn{2}{c}{4}        		& m    			\rule[1.2ex]{0pt}{1.2ex}\\
Pipe width                    	& $W$       			     	& \multicolumn{2}{c}{0.5}      		& m    		 	\rule[1.2ex]{0pt}{1.2ex}\\
Fluid density                  	& $\rho_{\rm f}$     	      	& \multicolumn{2}{c}{1.2}      		& $\rm kg/m^{3}$  	\rule[1.2ex]{0pt}{1.2ex} \\
Particle density             	& $\rho_{\rm p}$  			& \multicolumn{2}{c}{665, 2990}   	& $\rm kg/m^{3}$  	\rule[1.2ex]{0pt}{1.2ex}\\
Particle diameter           	& $d_{\rm p}$       			& \multicolumn{2}{c}{200, 400}    	& µm    		 	\rule[1.2ex]{0pt}{1.2ex}\\
Inlet particle volume fraction  	& $\phi^{\rm in}_{\rm p}$ 	& \multicolumn{2}{c}{0.005}     &     		-		 \rule[1.2ex]{0pt}{1.2ex}\\ \specialrule{.1em}{.05em}{.05em} 
\end{tabular}}
\end{table}
\begin{table}[H]\centering
\caption{Default physical parameters for the calculation of $\mathit{St_0}$ and the optimization problems.}
\label{opt1}
\fontsize{1}{1}\selectfont
\resizebox{1.0\linewidth}{!}{
\begin{tabular}{@{}lccll@{}}
																							 \specialrule{.1em}{.05em}{.05em} 
Parameters                   		& \multicolumn{1}{l}{Symbols} 	& \multicolumn{2}{c}{Values}   	& Unit 			\rule[1.2ex]{0pt}{1.2ex}\\ \specialrule{.1em}{.05em}{.05em} 
Inlet length               	    	& $L$       				& \multicolumn{2}{c}{0.33}      	& m    			\rule[1.2ex]{0pt}{1.2ex}\\
Maximum inlet fluid velocity          & $u^{\rm in}_{\rm f}$       	& \multicolumn{2}{c}{1.0}        	& m/s    			\rule[1.2ex]{0pt}{1.2ex}\\
Maximum inlet particle velocity     & $u^{\rm in}_{\rm p}$       	& \multicolumn{2}{c}{1.0}        	& m/s    			\rule[1.2ex]{0pt}{1.2ex}\\
Fluid density                  		& $\rho_{\rm f}$       		& \multicolumn{2}{c}{1.2}      	& $\rm kg/m^{3}$    	\rule[1.2ex]{0pt}{1.2ex}\\
Fluid viscosity                  	& $\mu$		       		& \multicolumn{2}{c}{0.04}      	& $\rm Pa\cdot s$    	\rule[1.2ex]{0pt}{1.2ex}\\
Particle density              		& $\rho_{\rm p}$       		& \multicolumn{2}{c}{1000} 	& $\rm kg/m^{3}$   	\rule[1.2ex]{0pt}{1.2ex}\\
Particle diameter           		& $d_{\rm p}$       		& \multicolumn{2}{c}{1}  	& mm    			\rule[1.2ex]{0pt}{1.2ex}\\
Inlet particle volume fraction  	& $\phi^{\rm in}_{\rm p}$       	& \multicolumn{2}{c}{0.005}    	& -    				\rule[1.2ex]{0pt}{1.2ex}\\ \specialrule{.1em}{.05em}{.05em} 
\end{tabular}}
\end{table}

\subsection{Particle drag variation maximization problem}
\subsubsection{Details of numerical setting}
In the flow simulation setup, the analysis domains of the two cases are discretized using the quadrilateral meshes of $80 \times 40$ and $80 \times 80$. In the optimization setting, the initial value of the design variables $\gamma$ is set to 0.5. The filter radius $R$ is set to 0.05. Additionally, the maximum value of the projection parameter $\beta_{ \rm max}$ and the projection threshold $\eta$ are set to 32 and 0.5, respectively. For the representation of the solid region, the maximum value of inverse permeability ${\overline{\alpha}}_{\rm f}$ is set to $1 \times 10^{4}$, and the convex parameter $q$ is set to 0.02. In the optimization problem, the weighting factors $w_{\Phi}$ and $w_{\zeta}$ in the objective function are set to 0.5, respectively, and the maximum limit of volume constraint $\overline{V}$ is set to $0.5V_{0}$. It should be noted that $c_{\alpha}$ is set to $\mathit{St}/\mathit{St}_{0}$, where $\mathit{St}$ represents the Stokes number, and $\mathit{St}_{0}$ corresponds to the Stokes number calculated according to Table ~\ref{opt1}. The expression of Stokes number is shown as follows:
\begin{equation}
\label{Steq}
\mathit{St} = \frac{\tau_{\rm p}}{\tau_{\rm f}} = \frac{\rho_{\rm p} d^{2}_{\rm p} U}{18 \mu_{\rm f} l_{\rm o}},
\end{equation}

where  $\tau_{\rm p}$ is the particle relaxation time, $\tau_{\rm f}$ is the characteristic time of the fluid flow, and $l_{\rm o}$ is the characteristic length of the obstacle. The Stokes number describes the behavior of particles suspended in a fluid flow. It represents the ratio of the particle relaxation time to the characteristic time of the fluid flow, indicating how easily particles can follow the flow.  As mentioned in section.~\ref{Fluid and solid representation using the design variable field}, as the inertia of particles increases, particles tend to follow their original trajectory rather than the fluid flow, which requires applying a larger fictitious body force to prevent particles from penetrating the solid region. Therefore, to properly account for the interaction between particles and the solid region, We set $c_{\alpha}$ based on the ratio of the Stokes number under different conditions to the reference Stokes number $\mathit{St}_{0}$.

\subsubsection{Case study 1: symmetrical flow setup}
\label{symmrtrical flow}
\begin{figure}[h]
\centering
\includegraphics[width=1.0\textwidth]{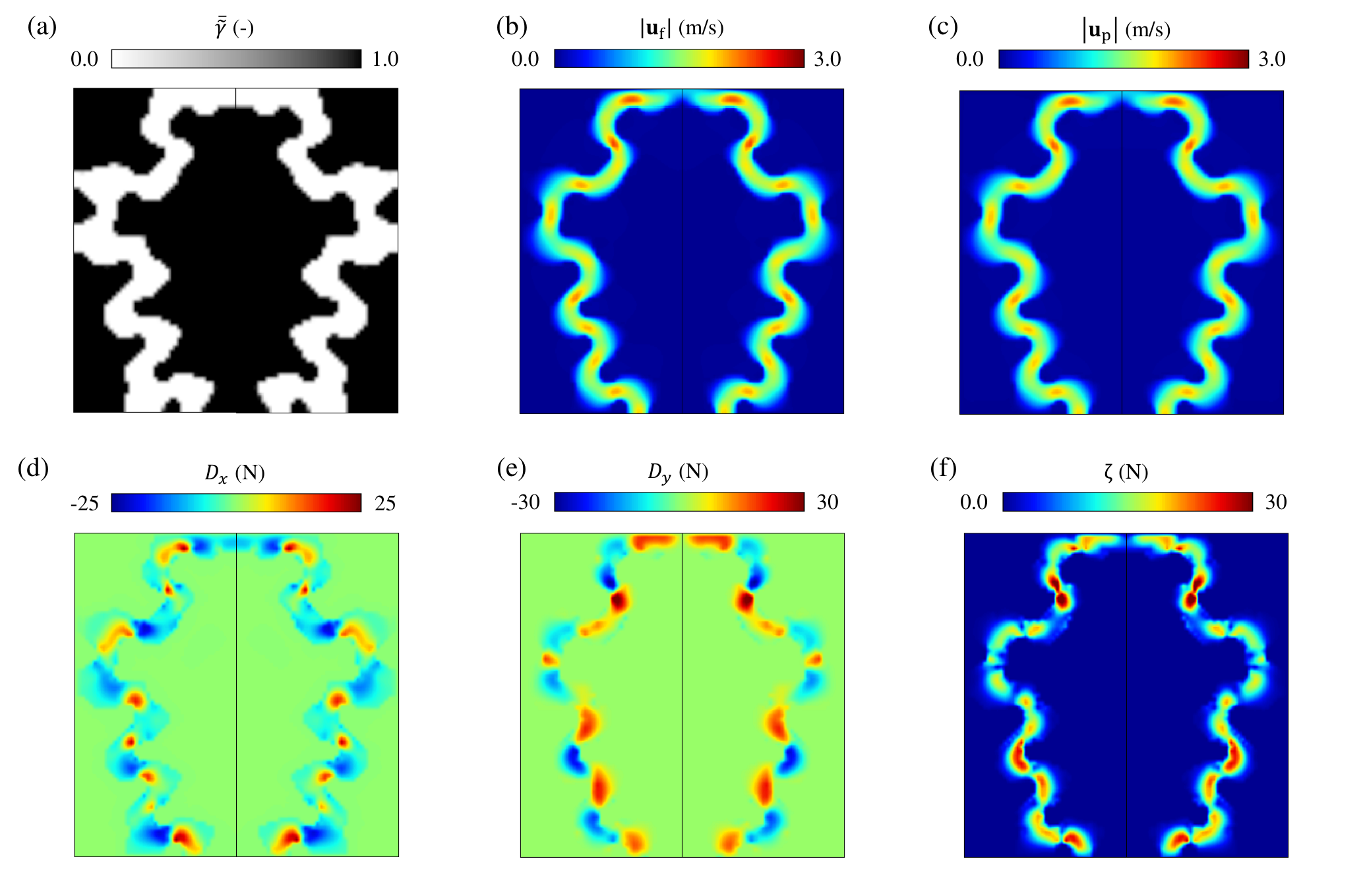} 
\caption{Distribution of design variables and physical quantities: (a) design variable $\bar{\tilde{\gamma}}$ (b) magnitude of fluid velocity $ \bold{u}_{\rm f} $ (c) magnitude of particle velocity $ \bold{u}_{\rm p} $ (d) magnitude of drag force in the $x$-direction $D_{x}$ (e) magnitude of drag force in the $y$-direction $D_{y}$ (f) magnitude of particle drag variation $\zeta$.}
\label{Fig.7}
\end{figure}

The optimized flow field is obtained in Fig.~\ref{Fig.7} according to the parameters in Table ~\ref{opt1}. The corresponding flow Reynolds number is 10. From Fig.~\ref{Fig.7}(a), it can be observed that the optimized flow field forms a serpentine flow field. Due to the relatively low Stokes number in this case, the distribution of the fluid velocity and particle velocity are nearly identical in Fig.~\ref{Fig.7}(b) and (c). However, in Fig.~\ref{Fig.7}(d) and (e), the drag force resulting from the velocity difference between the fluid and the particles is observed. Maximizing the variation in drag force on particles implies making the direction and magnitude of the particle velocity as different as possible from the fluid velocity. The curved structure of the serpentine flow channel induces a velocity difference between the fluid and the particles, which applies drag force on the particles. The distribution of drag force variation occurs at the turning points of the curved flow channels, as shown in Fig.~\ref{Fig.7}(f). 

More specifically, as particles enter a turning point in the $y$-direction, due to the obstruction of the curved structure, the velocity of the particles in the $y$-direction is lower than that of the fluid (indicated by the red region in Fig.~\ref{Fig.7}(e)), but higher in the $x$-direction (indicated by the blue region in Fig.~\ref{Fig.7}(d)); however, as particles leave the turning point, they continue to accelerate in the $y$-direction due to inertia and gravity effects, leading to their velocity exceeding that of the fluid in the $y$-direction (indicated by the blue region in Fig.~\ref{Fig.7}(e)). Meanwhile, their velocity in the $x$-direction will be lower than that of the fluid (indicated by the red region in Fig.~\ref{Fig.7}(d)). Through the repetitive bending of the serpentine flow path, particles will constantly accelerate and decelerate, causing the drag force on particles to change continuously. 

\begin{figure}[H]
\centering
\includegraphics[width=1.0\textwidth]{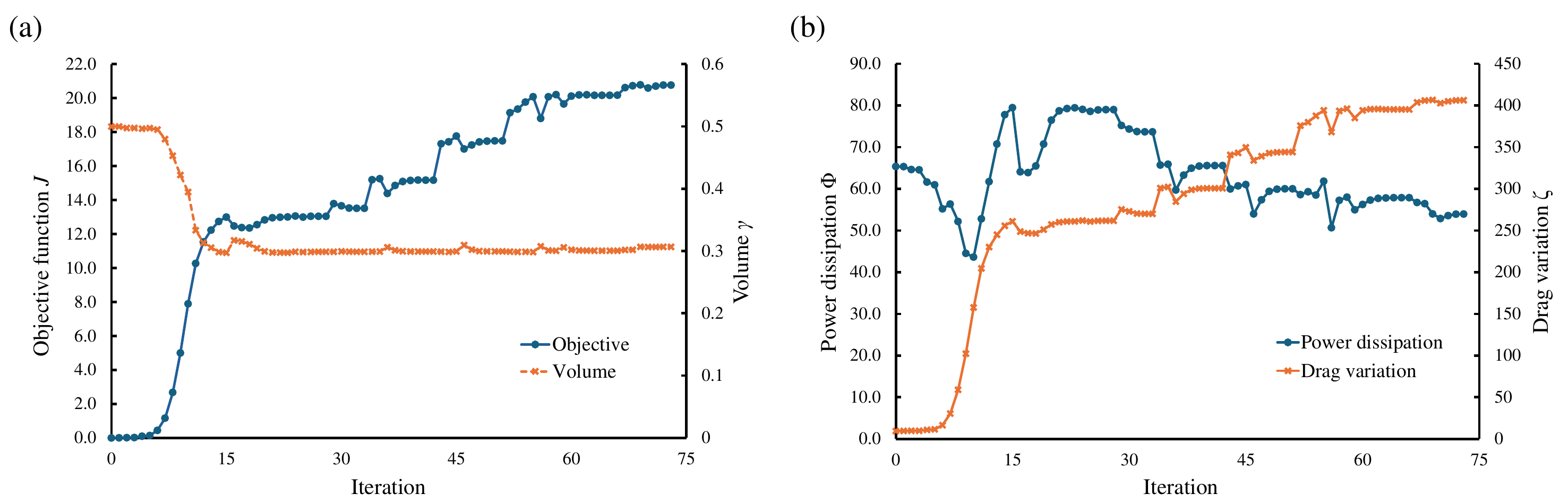} 
\caption{Convergence history of the optimized flow field in Fig.~\ref{Fig.7}: (a) objective and constraint (b) power dissipation and drag variation.}
\label{Fig.8}
\end{figure}

The convergence history is demonstrated in Fig.~\ref{Fig.8}. In addition to satisfying the objective function and constraints, the decrease in power dissipation and the increase in drag variation also indicate that the optimized flow field meets the requirement of increasing the resistance experienced by particles effectively.

\begin{figure}[H]
\centering
\includegraphics[width=0.9\textwidth]{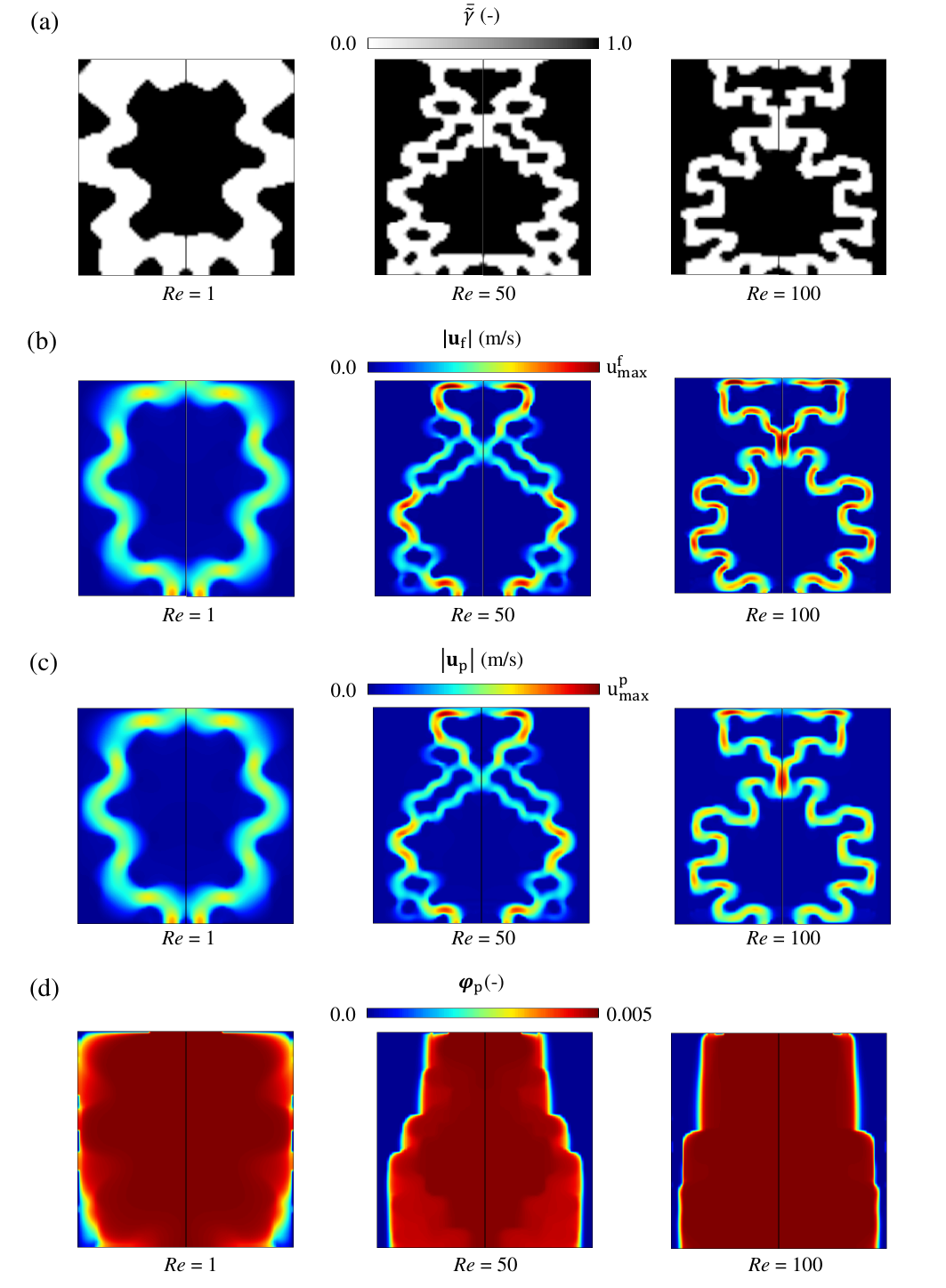} 
\caption{Optimized Results and physical quantities for different Reynolds numbers: (a) optimized flow field (b) fluid velocity (c) particle velocity (d) particle volume fraction. The maximum values on the color bar for fluid velocity and particle velocity are 2.0 $\rm m/s$, 3.0 $\rm m/s$, and 3.5 $\rm m/s$ with $Re$ = 1, 50, 100.}
\label{Fig.9}
\end{figure}

Furthermore, we examine the effect of different Reynolds numbers and Stokes numbers on the optimized flow field. In this investigation, the desired value of the Reynolds number is determined by adjusting the viscosity. The optimized flow fields with different Reynolds numbers are shown in Fig.~\ref{Fig.9}. As the Reynolds number increases, the optimized flow field becomes more complex, with an observable increase in the curvature of the flow channels. In the case of a smaller Reynolds number, since viscosity primarily dominates the flow behavior, increasing the drag force by bending the flow channels results in more energy loss. As a result, the number and curvature of bends in the serpentine flow fields are reduced. Moreover, due to the weak inertial effects on the particles, there is a slight difference in the velocity distribution between the particles and the fluid. In the case of a higher Reynolds number,  the inertial effects on particles become more significant, resulting in more noticeable differences in the velocity distribution between the particles and the fluid. Particles tend to collide with solids, changing the direction instead of following the fluid streamlines. This behavior increases the degree of variation in drag force.

It should be noted that due to the limitations of the fictitious body force formulation, the particle volume fraction shown in Fig.~\ref{Fig.9}(d) remains within the solid region. Nevertheless, since the fluid velocity and particle velocity inside the solid region are very low, the drag variation within the solid region is also minimal.

\begin{figure}[h]
\centering
\includegraphics[width=0.75\textwidth]{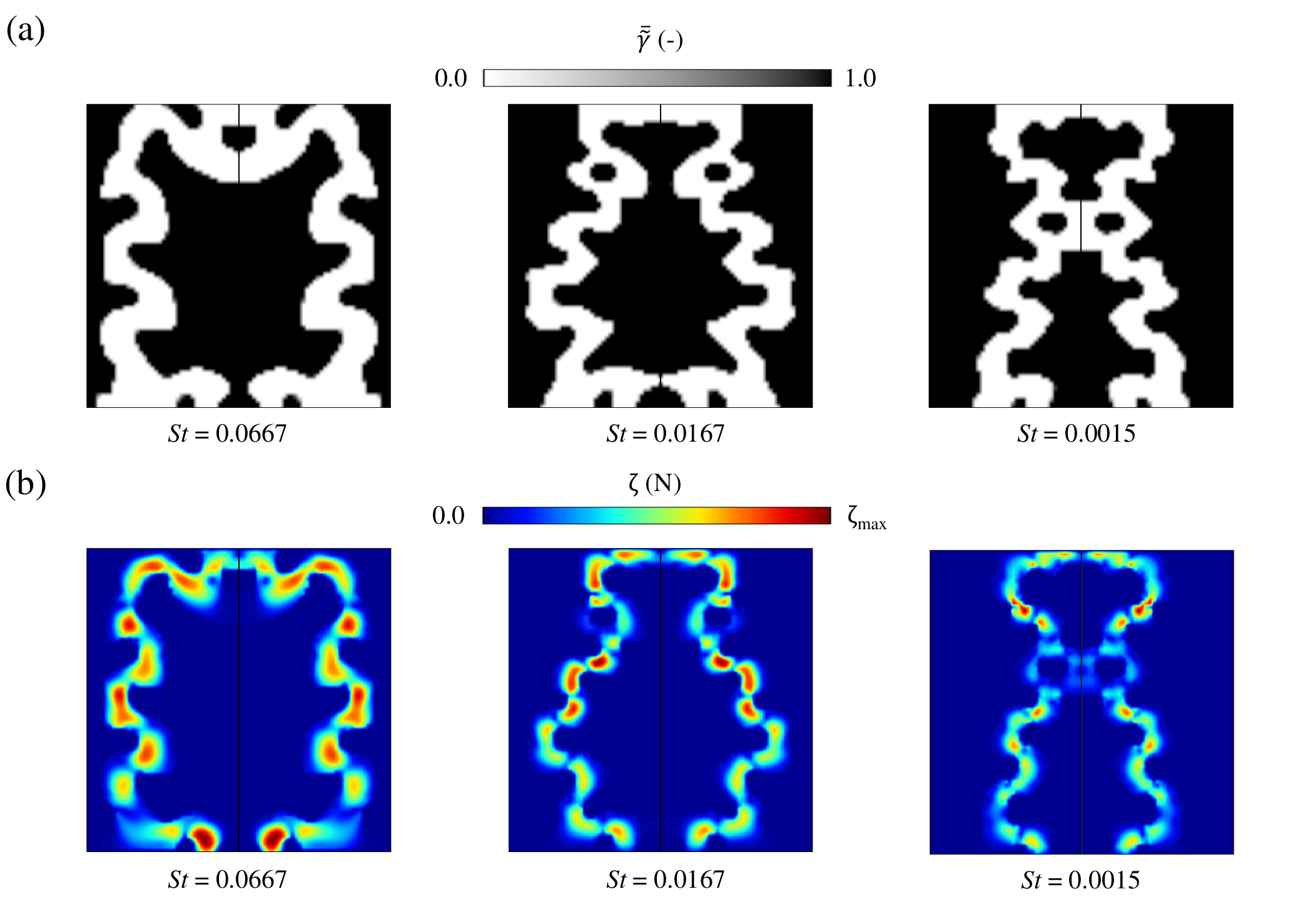} 
\caption{Optimized results and physical quantities for different Stokes numbers by adjusting the particle diameters to 4 mm, 2 mm, and 0.6 mm: (a) optimized flow field (b) magnitude of particle drag variation. The maximum values of the color bar for drag variation are 51.8, 27.2, and 3.77 with $\mathit{St} =$ 0.0667, 0.01667, and 0.0015, respectively.}
\label{Fig.10}
\end{figure}

\begin{figure}[H]
\centering
\includegraphics[width=0.75\textwidth]{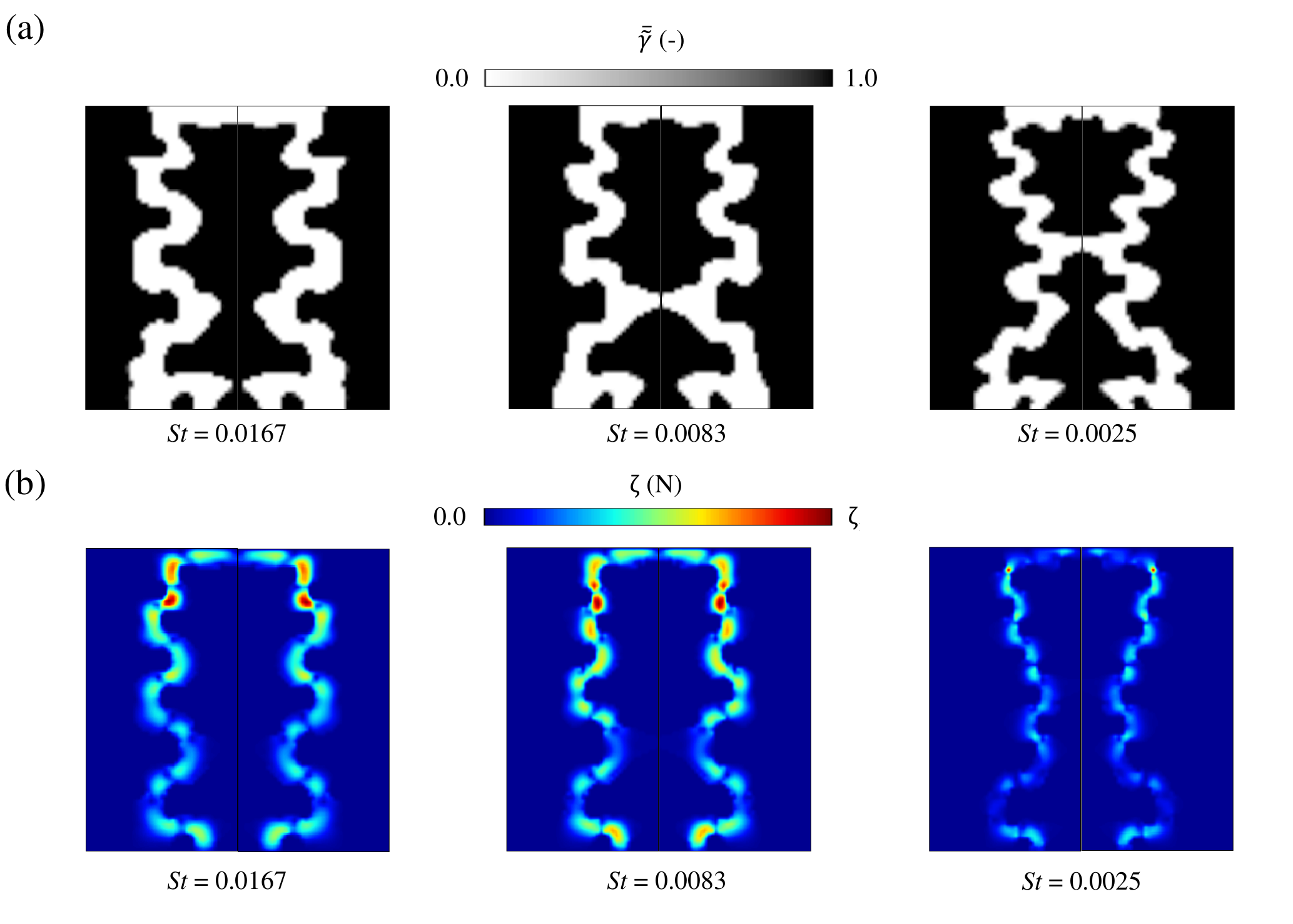} 
\caption{Optimized flow fields and physical quantities for different Stokes numbers by adjusting the particle densities $\rho_{\rm p}$ with 4000 $\rm kg/m^3$, 2000 $\rm kg/m^3$, and 600 $\rm kg/m^3$: (a) optimized flow field (b) the magnitude of particle drag variation. The maximum values of the color bar for drag variation are 37.75, 20.1, 16.7 with $\mathit{St}$ = 0.01667, 0.0083, and 0.0025, respectively.}
\label{Fig.11}
\end{figure}

Fig.~\ref{Fig.10} and Fig.~\ref{Fig.11} show the optimized configurations for different Stokes numbers, which correspond to various particle diameters and densities. Generally, the number of bends in the serpentine flow field increases as the Stokes number decreases. On the other hand, regarding bend curvature, a serpentine flow field can be regarded as a sequence of bent pipes, where the characteristic length of a bent pipe $l_{\rm o}$ is defined as its curvature radius $R_{\rm c}$~\cite{nicolaou2016characterization,henriquez2021numerical}. As the Stokes number decreases, the overall radius of curvature across all bends in the serpentine flow field decreases.

As the Stokes number decreases, particles follow the fluid streamlines more closely due to the shorter particle relaxation time. Consequently, the optimized flow field requires more bends and smaller curvature radii to enhance variations in the particle drag force. Conversely, as the Stokes number increases, the inertial effects on the particles are enhanced, causing them to adhere more closely to their original trajectories. It becomes easier for particles to detach from the flow streamlines with higher Stokes numbers. In such situations, the influence of bending structures on power dissipation exceeds their effect on particle drag force. The optimized flow field exhibits increased curvature radii and fewer bends to reduce power dissipation.

\subsubsection{Case study 2: asymmetrical flow setup}
\begin{figure}[h]
\centering
\includegraphics[width=0.85\textwidth]{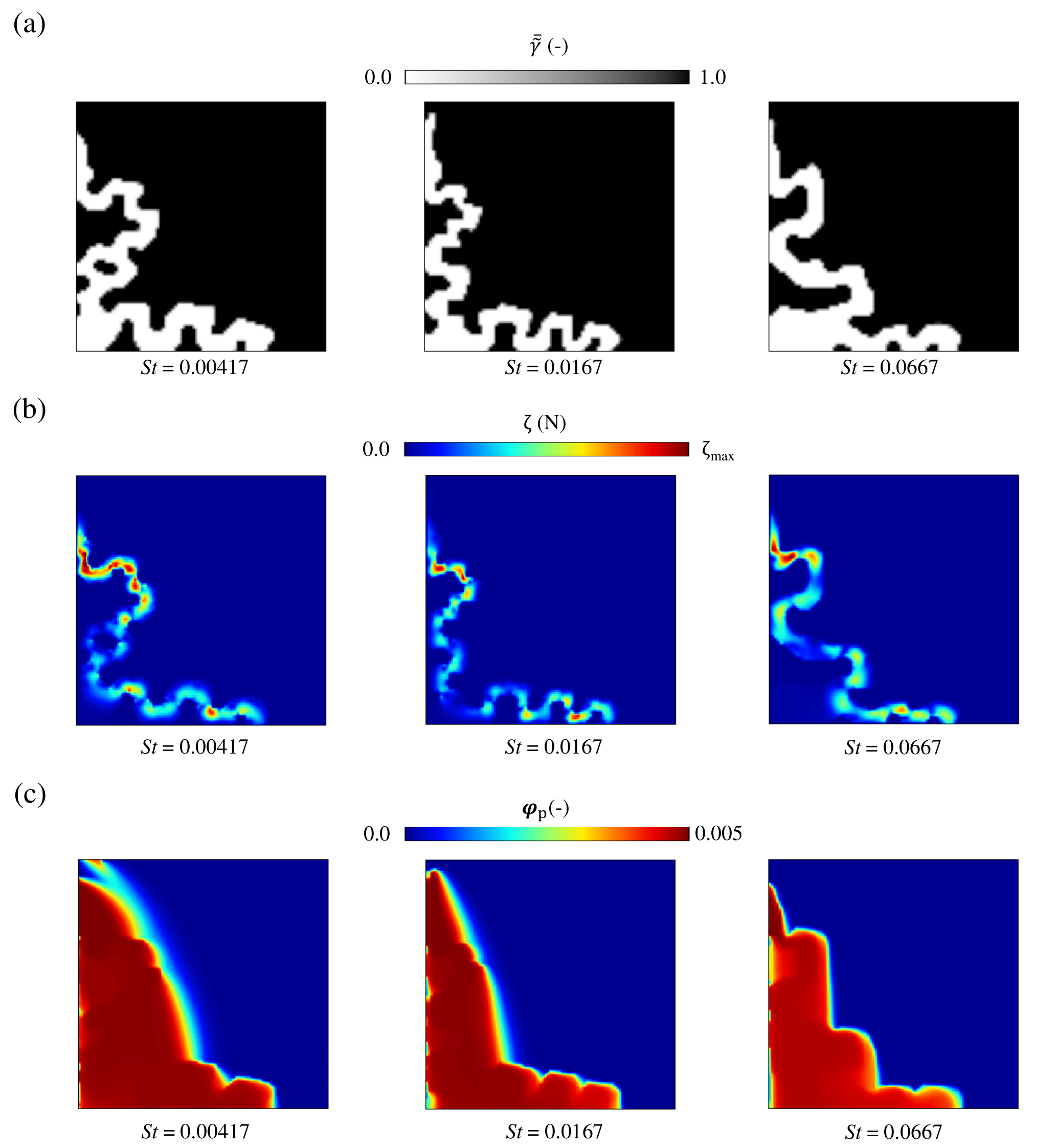} 
\caption{Optimized Results and physical quantities for different Stokes numbers by adjusting the particle diameters to 1 mm, 2 mm, and 4 mm: (a) optimized flow field (b) magnitude of particle drag variation (c) particle volume fraction. The maximum values of the color bar for drag variation are 20.4, 74.3, and 110.48 with $\mathit{St} =$ 0.00417, 0.01667, and 0.0667, respectively.}
\label{Fig.16}
\end{figure}

In this problem setup, we modify the layout of the inlet and outlet to investigate the effect of asymmetric flow on the optimized flow field. The optimization is also based on the parameters in Table ~\ref{opt1}. Fig.~\ref{Fig.16} shows the optimized configurations for different Stokes numbers by varying the particle diameter, and the serpentine feature of the flow path remains observable. As the Stokes number increases, both the number of bends and the overall radius of curvature across all bends in the flow field decrease, following the same trend observed in the first problem setting.

Notably, the overall trajectory of the serpentine flow paths in Fig.~\ref{Fig.16}(a) does not follow the shortest route. Instead, they first move downward vertically, then extend toward the outlet on the right. This L-shaped path suggests the presence of gravitational influence on particle movement.  Fig.~\ref{Fig.16}(c) shows that the particles tend to flow downward due to the effect of gravity, which explains the L-shaped trajectory of the serpentine flow path.

\begin{figure}[h]
\centering
\includegraphics[width=1.0\textwidth]{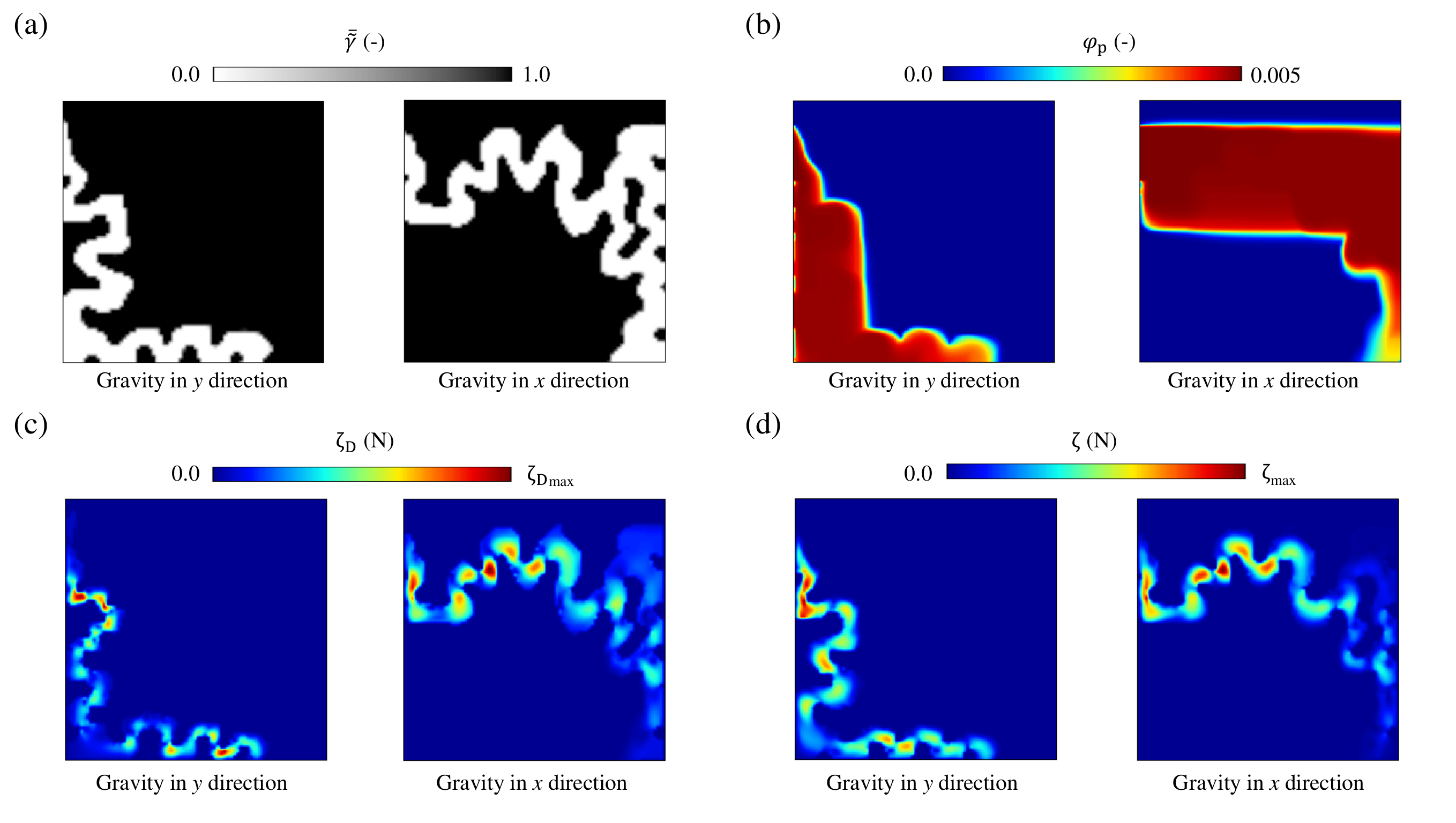} 
\caption{Optimized Results and physical quantities with $Re$ = 10, $d_{\rm p}$ = 2 mm, and $\rho_{\rm p}$ = 4000 $\rm kg/m^3$: (a) optimized flow field (b) particle volume fraction (c) magnitude of drag force (d) magnitude of particle drag variation. The maximum values of the color bar for drag force are 113.41 and 108.64 for gravity in $y$-direction and $x$-direction. The maximum values of the color bar for drag variation are 115.03 and 102.08 for gravity in $y$-direction and $x$-direction.}
\label{Fig.17}
\end{figure}

To further clarify the effect of gravity, Fig.~\ref{Fig.17} shows the optimized results under different gravity directions. As shown in Fig.~\ref{Fig.17}(a), although the overall trajectory of the serpentine flow path remains L-shaped under different gravity directions, the flow field direction is entirely different. Due to the effect of gravity, the particles move along the direction of gravity, and the resulting distribution of particle volume fraction significantly influences the optimized flow field, as illustrated in Fig.~\ref{Fig.17}(b). When particles are concentrated along the gravity direction, the particle volume fraction in other areas within the design domain becomes low, leading to a reduction in the effect of drag variation in those regions. Consequently, the overall trajectory of the optimized flow field aligns with the particle volume fraction distribution. In Fig.~\ref{Fig.17}(c), the drag distribution represents the velocity difference between the particles and the fluid, which primarily occurs along the serpentine flow path in the direction of gravity. Under the influence of gravity, particles tend to accelerate downward but are repeatedly decelerated as they are obstructed by the serpentine flow path. This continuous cycle of acceleration and deceleration increases the drag variation, as shown in Fig.~\ref{Fig.17}(d). Once the particles descend to the bottom, they can no longer accelerate under gravity, which reduces the velocity difference between the particles and the fluid, and consequently decreases drag variation. Nevertheless, due to the difference in relaxation times between the particles and the fluid, the optimized flow fields still present serpentine flow paths to enhance particle drag variation.

\section{Conclusion}
\label{Conclusion}
In this paper, we proposed a density-based topology optimization method for a particle flow problem with the Eulerian-Eulerian approach. The finite difference method named the NAPPLE algorithm was employed for fluid simulation on a collocated grid. We formulated the optimization problem as a multi-objective function, aiming to maximize the magnitude of particle drag variation while minimizing power dissipation. In the topology optimization framework, automatic differentiation was employed to compute discrete sensitivities without continuous adjoint formulation. The distribution of design variables was obtained via the GCMMA algorithm in the optimization process. For the numerical validation, some numerical examples of one-way coupling particle flow were conducted to compare with the prior works. For the numerical investigation, we considered two cases corresponding to symmetric and asymmetric flow and examined the effects of the Reynolds number, Stokes number, and gravity on the optimized flow field. Our key findings are summarized as follows:

\begin{itemize}
\item Flow fields with serpentine features effectively maximize the variation in particle drag force.
\item As the Reynolds number increases and the viscous effect decreases, the curvature of bends in the serpentine flow field increases. This enhancement raises the variation in particle drag force.
\item With the decrease in the Stokes number, due to the reduced particle relaxation time, the curvature and number of bends in the serpentine flow field increase to raise the variation in particle drag force.
\item As the gravitational effect increases, particles concentrate along the direction of gravity, leading the serpentine flow path to align with the particle distribution.
\end{itemize}

In future work, to extend this method to industrial applications, such as particle heating receivers or microfluidic devices, it is necessary to extend the numerical model to dense particle conditions and address the issue that the representation of fictitious forces does not fully account for the shear stress experienced by particles at the solid walls.
\section*{Acknowledgement}
This work was supported by JSPS KAKENHI (Grant No. 23H01323).

\appendix
\section{Validation of discrete sensitivities}
\label{sensitivity appendix}
\begin{figure}[h]
\centering
\includegraphics[width=0.5\textwidth]{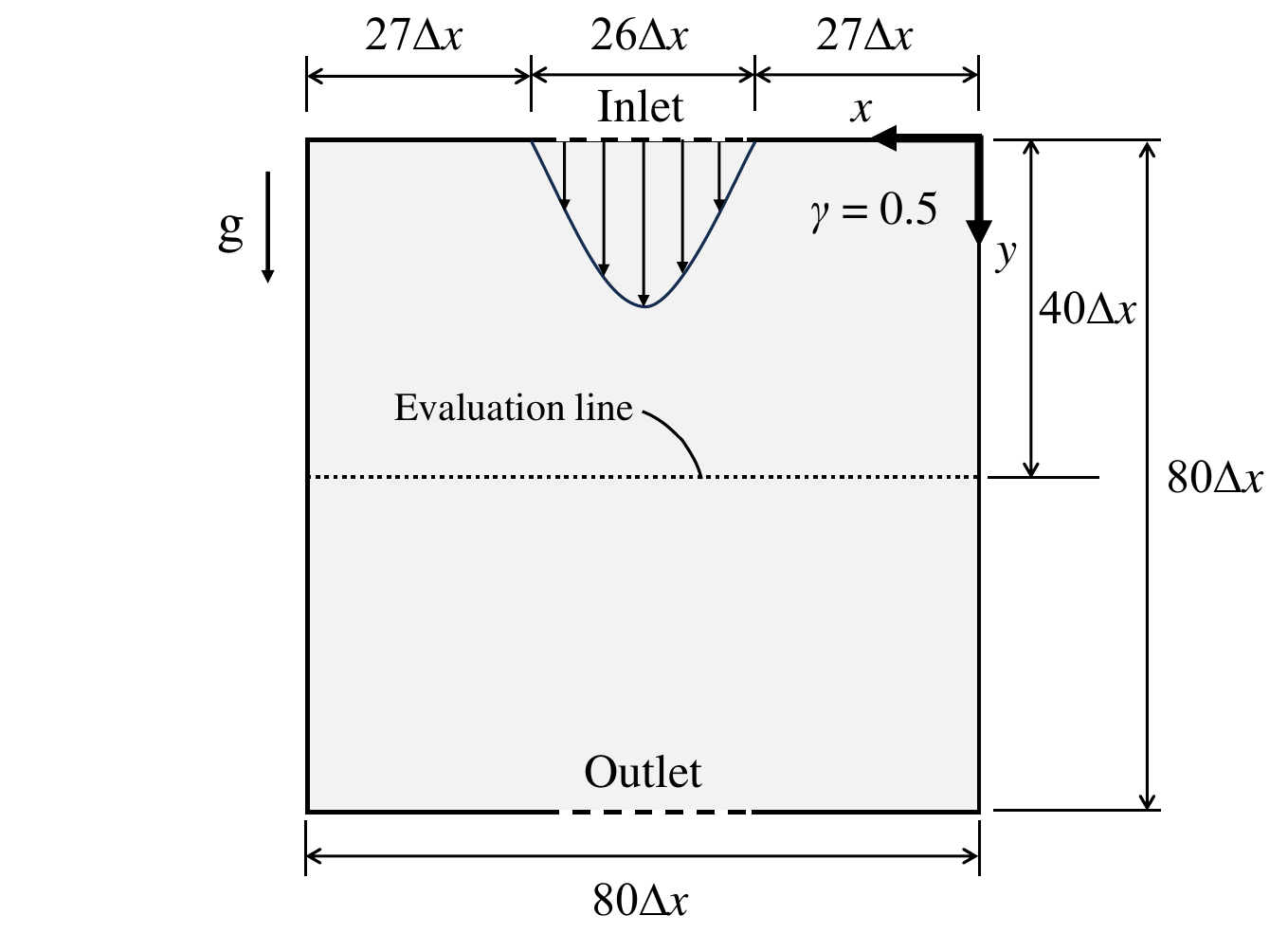} 
\caption{A schematic diagram of the evaluation line in the analysis domain to validate discrete sensitivities.}
\label{Fig.12}
\end{figure}

\begin{figure}[h]
\centering
\includegraphics[width=1.0\textwidth]{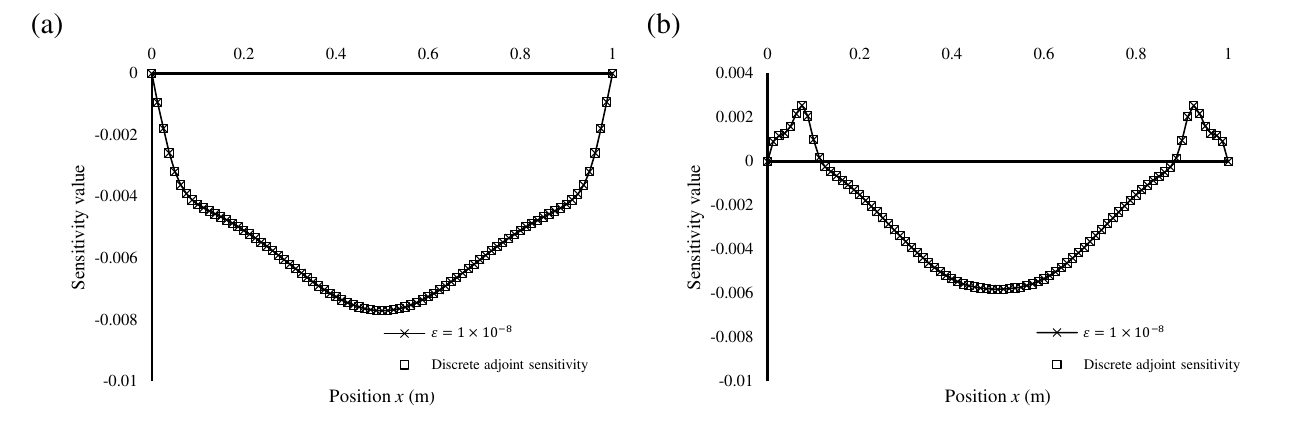} 
\caption{Comparison of sensitivities derived by automatic differentiation and finite difference check: (a) power dissipation $\Phi$ ($w_{\Phi} = 1, w_{\zeta} = 0$) (b) deviation of drag force $\zeta$ ($w_{\Phi} = 0, w_{\zeta} = 1$).}
\label{Fig.14}
\end{figure}

In this study, we implement automatic differentiation for the sensitivity calculation. To verify the correctness of the sensitivity calculation, we validate the sensitivities with the results of finite difference approximation as follows:

\begin{equation}
J^{'} = \frac{J(\gamma + \varepsilon) - J(\gamma - \varepsilon)}{2 \varepsilon},
\end{equation}
where $J^{'}$ is the sensitivities derived by finite difference approximation, and $\varepsilon$ is a small positive number. The sensitivities derived by automatic differentiation should match the outcome obtained by finite difference approximation. Fig.~\ref{Fig.12} illustrates the position of the evaluation line in the analysis domain. We examine the sensitivities for the power dissipation $\Phi$ and the variation of drag force $\zeta$ on the evaluation line separately, corresponding to two cases: $w_{\Phi} = 1, w_{\zeta} = 0$ and $w_{\Phi} = 0, w_{\zeta} = 1$. The physical parameters for validation are listed in Table~\ref{opt1}. As shown in Fig.~\ref{Fig.14}, the sensitivities obtained from automatic differentiation are entirely consistent with the results from the finite difference approximation. 

\begin{figure}[H]
\centering
\includegraphics[width=1.1\textwidth]{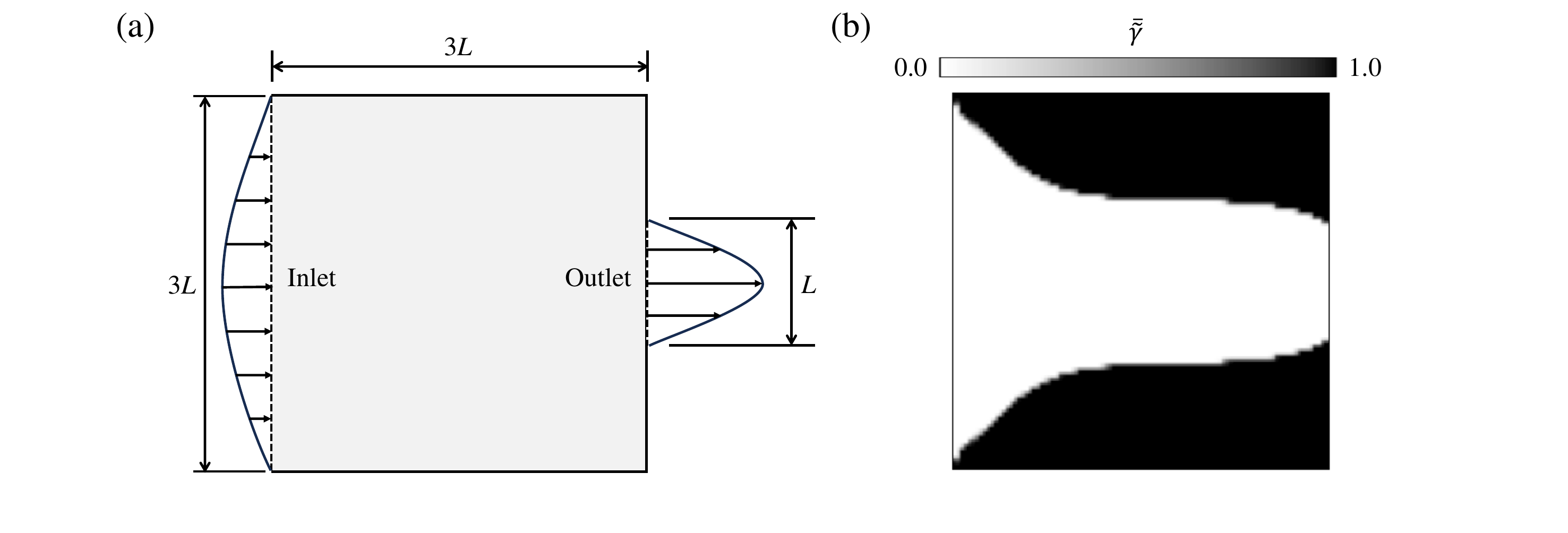} 
\caption{Analysis domain and optimized configuration for the diffuser problem: (a) analysis domain (b) optimized configuration.}
\label{Fig.13}
\end{figure}

In addition to the basic finite difference approximation, we further implement the diffuser problem to verify the accuracy of our sensitivity analysis. This problem is often regarded as a benchmark example in topology optimization studies~\cite{yaji2016topology, xie2021topology, suarez2018topology, padhy2023fluto}. Fig.~\ref{Fig.13}(a) shows the analysis domain. In this optimization problem, we employ the power dissipation $\Phi$ for the objective function and set the maximum value of volume constraint to $0.5V_{0}$ ($V_{0}$ is the volume of the design domain $D$). The flow Reynolds number is set to 1.0 in the numerical simulation. The velocity boundary conditions for both the inlet and the outlet are parabolic profiles, with the maximum magnitudes set to 1.0 at the inlet and 3.0 at the outlet. Fig.~\ref{Fig.13}(b) shows the optimized configuration derived through this method, which matches the prior research~\cite{borrvall2003topology}.


\begin{thebibliography}{10}
\expandafter\ifx\csname url\endcsname\relax
  \def\url#1{\texttt{#1}}\fi
\expandafter\ifx\csname urlprefix\endcsname\relax\def\urlprefix{URL }\fi
\expandafter\ifx\csname href\endcsname\relax
  \def\href#1#2{#2} \def\path#1{#1}\fi

\bibitem{klaewkla2011review}
R.~Klaewkla, M.~Arend, W.~F. Hoelderich, A review of mass transfer controlling
  the reaction rate in heterogeneous catalytic systems, Vol.~5, INTECH Open
  Access Publisher Rijeka, 2011.

\bibitem{duan2015green}
H.~Duan, D.~Wang, Y.~Li, Green chemistry for nanoparticle synthesis, Chemical
  Society Reviews 44~(16) (2015) 5778--5792.

\bibitem{duduta2011semi}
M.~Duduta, B.~Ho, V.~C. Wood, P.~Limthongkul, V.~E. Brunini, W.~C. Carter,
  Y.-M. Chiang, Semi-solid lithium rechargeable flow battery, Advanced Energy
  Materials 1~(4) (2011) 511--516.

\bibitem{wang2020review}
H.~Wang, J.~Lu, A review on particle size effect in metal-catalyzed
  heterogeneous reactions, Chinese Journal of Chemistry 38~(11) (2020)
  1422--1444.

\bibitem{santana2019development}
H.~S. Santana, J.~L. Silva~Jr, O.~P. Taranto, Development of microreactors
  applied on biodiesel synthesis: {From} experimental investigation to
  numerical approaches, Journal of industrial and engineering chemistry 69
  (2019) 1--12.

\bibitem{suh2010review}
Y.~K. Suh, S.~Kang, A review on mixing in microfluidics, Micromachines 1~(3)
  (2010) 82--111.

\bibitem{chaurasia2020thermal}
S.~R. Chaurasia, R.~Sarviya, Thermal performance analysis of {CuO/water}
  nanofluid flow in a pipe with single and double strip helical screw tape,
  Applied Thermal Engineering 166 (2020) 114631.

\bibitem{das2006heat}
S.~K. Das, S.~U. Choi, H.~E. Patel, Heat transfer in nanofluids—a review,
  Heat transfer engineering 27~(10) (2006) 3--19.

\bibitem{hong2017millifluidic}
L.~Hong, T.-L. Cheung, N.~Rao, Q.~Ouyang, Y.~Wang, S.~Zeng, C.~Yang, D.~Cuong,
  P.~H.~J. Chong, L.~Liu, et~al., Millifluidic synthesis of cadmium sulfide
  nanoparticles and their application in bioimaging, RSC advances 7~(58) (2017)
  36819--36832.

\bibitem{kehlenbeck2002particle}
R.~Kehlenbeck, J.~G. Yates, R.~Di~Felice, H.~Hofbauer, R.~Rauch, Particle
  residence time and particle mixing in a scaled internal circulating fluidized
  bed, Industrial \& engineering chemistry research 41~(11) (2002) 2637--2645.

\bibitem{phirommark2023cfd}
P.~Phirommark, S.~Namchanthra, J.~Chaiyanupong, S.~Uapipatanakul, W.~Chookaew,
  C.~Suvanjumrat, M.~Promtong, {CFD} elucidation of microscopic particles in a
  low-volumetric classifier towards effects of {Stokes} number and density
  ratio, International Journal of Thermofluids 20 (2023) 100497.

\bibitem{alhamdan1997residence}
A.~Alhamdan, S.~Sastry, Residence time distribution of food and simulated
  particles in a holding tube, Journal of food engineering 34~(3) (1997)
  271--292.

\bibitem{lee2018parametric}
T.~Lee, S.~Shin, S.~I. Abdel-Khalik, Parametric investigation of particulate
  flow in interconnected porous media for central particle-heating receiver,
  Journal of Mechanical Science and Technology 32 (2018) 1181--1186.

\bibitem{afrin2014modeling}
S.~Afrin, J.~D. Ortega, C.~K. Ho, V.~Kumar, {Modeling of a
  High-Temperature-Serpentine External Tubular Receiver Using Supercritical
  CO2}, in: Energy Sustainability, Vol. 45868, American Society of Mechanical
  Engineers, 2014, p. V001T02A011.

\bibitem{al2020study}
H.~Al-Ansary, A.~El-Leathy, A.~Alswaiyd, S.~Alaqel, N.~Saleh, R.~Saeed,
  Z.~Al-Suhaibani, S.~Danish, E.~Djajadiwinata, S.~Jeter, Study of the optimum
  discrete structure configuration in obstructed flow particle heating
  receivers, in: AIP Conference Proceedings, Vol. 2303, AIP Publishing, 2020.

\bibitem{bovskovic2011residence}
D.~Bo{\v{s}}kovi{\'c}, S.~Loebbecke, G.~Gross, J.~Koehler, Residence time
  distribution studies in microfluidic mixing structures, Chemical engineering
  \& technology 34~(3) (2011) 361--370.

\bibitem{yang2015experimental}
S.~Yang, H.~Li, Q.~Zhu, Experimental study and numerical simulation of baffled
  bubbling fluidized beds with {Geldart A} particles in three dimensions,
  Chemical Engineering Journal 259 (2015) 338--347.

\bibitem{kreesaeng2024effect}
S.~Kreesaeng, B.~Chalermsinsuwan, P.~Piumsomboon, {Effect of Inserting Baffles
  on the Solid Particle Segregation Behavior in Fluidized Bed Reactor: A
  Computational Study}, ChemEngineering 8~(1) (2024) 7.

\bibitem{rennebaum2024effect}
H.~S. Rennebaum, D.~L. Brummerloh, S.~Benders, A.~Penn, The effect of baffles
  on the hydrodynamics of a gas-solid fluidized bed studied using real-time
  magnetic resonance imaging, Powder technology 432 (2024) 119114.

\bibitem{an2012computational}
H.~An, A.~Li, A.~P. Sasmito, J.~C. Kurnia, S.~V. Jangam, A.~S. Mujumdar,
  Computational fluid dynamics {(CFD)} analysis of micro-reactor performance:
  {Effect} of various configurations, Chemical engineering science 75 (2012)
  85--95.

\bibitem{coyle2008micro}
E.~E. Coyle, M.~Oelgem{\"o}ller, Micro-photochemistry: photochemistry in
  microstructured reactors. {The} new photochemistry of the future?,
  Photochemical \& Photobiological Sciences 7 (2008) 1313--1322.

\bibitem{peres2019analysis}
J.~C.~G. Peres, C.~d.~C. Herrera, S.~L. Baldochi, W.~de~Rossi, A.~dos Santos
  Vianna~Jr, Analysis of a microreactor for synthesizing nanocrystals by
  computational fluid dynamics, The Canadian Journal of Chemical Engineering
  97~(2) (2019) 594--603.

\bibitem{ho2016review}
C.~K. Ho, A review of high-temperature particle receivers for concentrating
  solar power, Applied Thermal Engineering 109 (2016) 958--969.

\bibitem{ho2019sun}
C.~K. Ho, J.~M. Christian, J.~E.~Yellowhair, K.~Armijo, W.~J. Kolb, S.~Jeter,
  M.~Golob, C.~Nguyen, On-sun performance evaluation of alternative
  high-temperature falling particle receiver designs, Journal of Solar Energy
  Engineering 141~(1) (2019) 011009.

\bibitem{aubin2009effect}
J.~Aubin, L.~Prat, C.~Xuereb, C.~Gourdon, Effect of microchannel aspect ratio
  on residence time distributions and the axial dispersion coefficient,
  Chemical Engineering and Processing: Process Intensification 48~(1) (2009)
  554--559.

\bibitem{bendsoe1988generating}
M.~P. Bends{\o}e, N.~Kikuchi, Generating optimal topologies in structural
  design using a homogenization method, Computer methods in applied mechanics
  and engineering 71~(2) (1988) 197--224.

\bibitem{bendsoe2008topology}
M.~P. Bends{\o}e, O.~Sigmund, Topology optimization, Theory, Methods and
  (2008).

\bibitem{borrvall2003topology}
T.~Borrvall, J.~Petersson, {Topology optimization of fluids in Stokes flow},
  International journal for numerical methods in fluids 41~(1) (2003) 77--107.

\bibitem{guest2006topology}
J.~K. Guest, J.~H. Pr{\'e}vost, Topology optimization of creeping fluid flows
  using a {Darcy--Stokes} finite element, International Journal for Numerical
  Methods in Engineering 66~(3) (2006) 461--484.

\bibitem{gersborg2005topology}
A.~Gersborg-Hansen, O.~Sigmund, R.~B. Haber, Topology optimization of channel
  flow problems, Structural and multidisciplinary optimization 30 (2005)
  181--192.

\bibitem{kreissl2011topology}
S.~Kreissl, G.~Pingen, K.~Maute, Topology optimization for unsteady flow,
  International Journal for Numerical Methods in Engineering 87~(13) (2011)
  1229--1253.

\bibitem{deng2011topology}
Y.~Deng, Z.~Liu, P.~Zhang, Y.~Liu, Y.~Wu, Topology optimization of unsteady
  incompressible {Navier--Stokes} flows, Journal of Computational Physics
  230~(17) (2011) 6688--6708.

\bibitem{yaji2015topology}
K.~Yaji, T.~Yamada, S.~Kubo, K.~Izui, S.~Nishiwaki, A topology optimization
  method for a coupled thermal--fluid problem using level set boundary
  expressions, International Journal of Heat and Mass Transfer 81 (2015)
  878--888.

\bibitem{yaji2018topology}
K.~Yaji, S.~Yamasaki, S.~Tsushima, T.~Suzuki, K.~Fujita, Topology optimization
  for the design of flow fields in a redox flow battery, Structural and
  multidisciplinary optimization 57 (2018) 535--546.

\bibitem{chen2019computational}
C.~H. Chen, K.~Yaji, S.~Yamasaki, S.~Tsushima, K.~Fujita, Computational design
  of flow fields for vanadium redox flow batteries via topology optimization,
  Journal of Energy Storage 26 (2019) 100990.

\bibitem{jenkins2015level}
N.~Jenkins, K.~Maute, Level set topology optimization of stationary
  fluid-structure interaction problems, Structural and Multidisciplinary
  Optimization 52 (2015) 179--195.

\bibitem{li2023topology}
P.~Li, L.~Shi, J.~Zhao, B.~Liu, H.~Yan, Y.~Deng, B.~Yin, T.~Zhou, Y.~Zhu,
  Topology optimization design of a passive two-dimensional micromixer,
  Chemical Physics Letters 821 (2023) 140445.

\bibitem{andreasen2020framework}
C.~S. Andreasen, A framework for topology optimization of inertial microfluidic
  particle manipulators, Structural and Multidisciplinary Optimization 61~(6)
  (2020) 2481--2499.

\bibitem{yoon2022transient}
G.~H. Yoon, Transient sensitivity analysis and topology optimization of
  particle suspended in transient laminar fluid, Computer Methods in Applied
  Mechanics and Engineering 393 (2022) 114696.

\bibitem{yoon2021development}
G.~H. Yoon, H.~So, Development of topological optimization schemes controlling
  the trajectories of multiple particles in fluid, Structural and
  Multidisciplinary Optimization 63 (2021) 2355--2373.

\bibitem{choi2023matlab}
Y.~H. Choi, G.~H. Yoon, A {MATLAB topology} optimization code to control the
  trajectory of particle in fluid, Structural and Multidisciplinary
  Optimization 66~(4) (2023) 91.

\bibitem{zahari2018introduction}
N.~Zahari, M.~Zawawi, L.~M. Sidek, D.~Mohamad, Z.~Itam, M.~Ramli, A.~Syamsir,
  A.~Abas, M.~Rashid, Introduction of discrete phase model {(DPM)} in fluid
  flow: a review, in: AIP Conference Proceedings, Vol. 2030, AIP Publishing,
  2018.

\bibitem{balakin2010eulerian}
B.~V. Balakin, A.~C. Hoffmann, P.~Kosinski, L.~D. Rhyne, {Eulerian-Eulerian
  CFD} model for the sedimentation of spherical particles in suspension with
  high particle concentrations, Engineering Applications of Computational Fluid
  Mechanics 4~(1) (2010) 116--126.

\bibitem{zhang2019simulation}
Y.~Zhang, Z.~Ran, B.~Jin, Y.~Zhang, C.~Zhou, F.~Sher, Simulation of particle
  mixing and separation in multi-component fluidized bed using
  {Eulerian-Eulerian method: A review}, International Journal of Chemical
  Reactor Engineering 17~(11) (2019) 20190064.

\bibitem{gao2012review}
Y.~Gao, F.~J. Muzzio, M.~G. Ierapetritou, A review of the {Residence Time
  Distribution (RTD)} applications in solid unit operations, Powder technology
  228 (2012) 416--423.

\bibitem{ariyaratne2016cfd}
W.~H. Ariyaratne, E.~Manjula, C.~Ratnayake, M.~C. Melaaen, {CFD approaches for
  modeling gas-solids multiphase flows--A review}, in: 9th EUROSIM \& the 57th
  SIMS Conference, Oulu, 2016, pp. 680--686.

\bibitem{alexandersen2020review}
J.~Alexandersen, C.~S. Andreasen, A review of topology optimisation for
  fluid-based problems, Fluids 5~(1) (2020) 29.

\bibitem{yan2023topology}
K.~Yan, Y.~Wang, Y.~Pan, G.~Sun, J.~Chen, X.~Cai, G.~Cheng, Topology
  optimization of simplified convective heat transfer problems using the finite
  volume method, Science China Technological Sciences 66~(5) (2023) 1352--1364.

\bibitem{makhija2012topology}
D.~Makhija, G.~Pingen, R.~Yang, K.~Maute, Topology optimization of
  multi-component flows using a multi-relaxation time lattice {Boltzmann}
  method, Computers \& fluids 67 (2012) 104--114.

\bibitem{sasaki2019topology}
Y.~Sasaki, Y.~Sato, T.~Yamada, K.~Izui, S.~Nishiwaki, Topology optimization for
  fluid flows using the {MPS} method incorporating the level set method,
  Computers \& Fluids 188 (2019) 86--101.

\bibitem{zawawi2018review}
M.~H. Zawawi, A.~Saleha, A.~Salwa, N.~Hassan, N.~M. Zahari, M.~Z. Ramli, Z.~C.
  Muda, A review: Fundamentals of computational fluid dynamics {(CFD)}, in: AIP
  conference proceedings, Vol. 2030, AIP Publishing, 2018.

\bibitem{fanxi2017fast}
L.~Fanxi, X.~Tianhang, Y.~Xiongqing, A fast and automatic full-potential finite
  volume solver on {Cartesian} grids for unconventional configurations, Chinese
  Journal of Aeronautics 30~(3) (2017) 951--963.

\bibitem{barman2016introduction}
P.~C. Barman, Introduction to computational fluid dynamics, International
  Journal of Information Science and Computing 3~(2) (2016) 117--120.

\bibitem{anderson1995computational}
J.~D. Anderson, J.~Wendt, Computational fluid dynamics, Vol. 206, Springer,
  1995.

\bibitem{moukalled2016finite}
F.~Moukalled, L.~Mangani, M.~Darwish, F.~Moukalled, L.~Mangani, M.~Darwish, The
  finite volume method, Springer, 2016.

\bibitem{cheng2009high}
J.~Cheng, C.-W. Shu, {High order schemes for CFD: A review}, Chinese Journal of
  Computational Physics 26~(5) (2009) 633.

\bibitem{rizzetta2008high}
D.~P. Rizzetta, M.~R. Visbal, P.~E. Morgan, A high-order compact
  finite-difference scheme for large-eddy simulation of active flow control,
  Progress in Aerospace Sciences 44~(6) (2008) 397--426.

\bibitem{lee1992artificial}
S.~Lee, R.~Tzong, {Artificial pressure for pressure-linked equation},
  International journal of heat and mass transfer 35~(10) (1992) 2705--2716.

\bibitem{CCSA}
K.~Svanberg, A class of globally convergent optimization methods based on
  conservative convex separable approximations, {SIAM} Journal on Optimization
  12 (2002) 555--573.

\bibitem{griewank2008evaluating}
A.~Griewank, A.~Walther, Evaluating derivatives: principles and techniques of
  algorithmic differentiation, SIAM, 2008.

\bibitem{griewank2000algorithm}
A.~Griewank, A.~Walther, Algorithm 799: revolve: an implementation of
  checkpointing for the reverse or adjoint mode of computational
  differentiation, ACM Transactions on Mathematical Software (TOMS) 26~(1)
  (2000) 19--45.

\bibitem{crowe2011multiphase}
C.~T. Crowe, J.~D. Schwarzkopf, M.~Sommerfeld, Y.~Tsuji, {Multiphase flows with
  droplets and particles}, CRC press, 2011.

\bibitem{kim2010new}
M.-C. Kim, C.~Klapperich, A new method for simulating the motion of individual
  ellipsoidal bacteria in microfluidic devices, Lab on a Chip 10~(18) (2010)
  2464--2471.

\bibitem{wang2019industry}
Z.~Wang, C.~Liu, W.~Wei, Industry applications of magnetic separation based on
  nanoparticles: {A} review, International Journal of Applied Electromagnetics
  and Mechanics 60~(2) (2019) 281--297.

\bibitem{khashan2014computational}
S.~Khashan, A.~Alazzam, E.~Furlani, Computational analysis of enhanced magnetic
  bioseparation in microfluidic systems with flow-invasive magnetic elements,
  Scientific reports 4~(1) (2014) 5299.

\bibitem{hicdurmaz2022numerical}
S.~Hicdurmaz, E.~F. Johnson, J.~Grobbel, L.~Amsbeck, R.~Buck, B.~Hoffschmidt,
  Numerical heat transfer modelling of a centrifugal solar particle receiver,
  in: AIP Conference Proceedings, Vol. 2445, AIP Publishing, 2022.

\bibitem{sajeesh2014particle}
P.~Sajeesh, A.~K. Sen, Particle separation and sorting in microfluidic devices:
  a review, Microfluidics and nanofluidics 17 (2014) 1--52.

\bibitem{kotoky2018development}
S.~Kotoky, Development and application of a generic finite volume multiphase
  flow solver for gas particulate flows.

\bibitem{gidaspow1994multiphase}
D.~Gidaspow, Multiphase flow and fluidization: continuum and kinetic theory
  descriptions, Academic press, 1994.

\bibitem{bourdin2001filters}
B.~Bourdin, Filters in topology optimization, International journal for
  numerical methods in engineering 50~(9) (2001) 2143--2158.

\bibitem{bruns2001topology}
T.~E. Bruns, D.~A. Tortorelli, Topology optimization of non-linear elastic
  structures and compliant mechanisms, Computer methods in applied mechanics
  and engineering 190~(26-27) (2001) 3443--3459.

\bibitem{diaz1995checkerboard}
A.~Diaz, O.~Sigmund, Checkerboard patterns in layout optimization, Structural
  optimization 10 (1995) 40--45.

\bibitem{wang2011projection}
F.~Wang, B.~S. Lazarov, O.~Sigmund, On projection methods, convergence and
  robust formulations in topology optimization, Structural and
  multidisciplinary optimization 43 (2011) 767--784.

\bibitem{prakash1985control}
C.~Prakash, S.~Patankar, A control volume-based finite-element method for
  solving the {Navier-Stokes} equations using equal-order velocity-pressure
  interpolation, Numerical heat transfer 8~(3) (1985) 259--280.

\bibitem{date1996complete}
A.~Date, Complete pressure correction algorithm for solution of incompressible
  {Navier-Stokes} equations on a nonstaggered grid, Numerical Heat Transfer
  29~(4) (1996) 441--458.

\bibitem{ferziger2019computational}
J.~H. Ferziger, M.~Peri{\'c}, R.~L. Street, Computational methods for fluid
  dynamics, springer, 2019.

\bibitem{shong1989weighting}
L.~Shong-Leih, {Weighting function scheme and its application on
  multidimensional conservation equations}, International journal of heat and
  mass transfer 32~(11) (1989) 2065--2073.

\bibitem{patankar2018numerical}
S.~Patankar, Numerical heat transfer and fluid flow, Taylor \& Francis, 2018.

\bibitem{lee1988new}
S.-L. Lee, A new numerical formulation for parabolic differential equations
  under the consideration of large time steps, International journal for
  numerical methods in engineering 26~(7) (1988) 1541--1549.

\bibitem{patankar1983calculation}
S.~V. Patankar, D.~B. Spalding, A calculation procedure for heat, mass and
  momentum transfer in three-dimensional parabolic flows, in: Numerical
  prediction of flow, heat transfer, turbulence and combustion, Elsevier, 1983,
  pp. 54--73.

\bibitem{moukalled2003pressure}
F.~Moukalled, M.~Darwish, B.~Sekar, A pressure-based algorithm for multi-phase
  flow at all speeds, Journal of Computational Physics 190~(2) (2003) 550--571.

\bibitem{lee1990strongly}
S.-L. Lee, A strongly implicit solver for two-dimensional elliptic differential
  equations, Numerical Heat Transfer, Part B Fundamentals 16~(2) (1990)
  161--178.

\bibitem{hogan2014fast}
R.~J. Hogan, Fast reverse-mode automatic differentiation using expression
  templates in {C++}, ACM Transactions on Mathematical Software (TOMS) 40~(4)
  (2014) 1--16.

\bibitem{towara2013discrete}
M.~Towara, U.~Naumann, A discrete adjoint model for {OpenFOAM}, Procedia
  Computer Science 18 (2013) 429--438.

\bibitem{kotoky2018parametric}
S.~Kotoky, A.~Dalal, G.~Natarajan, A parametric study of dispersed laminar
  gas-particle flows through vertical and horizontal channels, Advanced Powder
  Technology 29~(5) (2018) 1072--1084.

\bibitem{nicolaou2016characterization}
L.~Nicolaou, T.~Zaki, Characterization of aerosol {Stokes} number in 90° bends
  and idealized extrathoracic airways, Journal of Aerosol Science 102 (2016)
  105--127.

\bibitem{henriquez2021numerical}
S.~Henr{\'\i}quez~Lira, M.~J. Torres, R.~Guerra~Silva, J.~Zahr~Vi{\~n}uela,
  Numerical characterization of the solid particle accumulation in a turbulent
  flow through curved pipes by means of stokes numbers, Applied Sciences
  11~(16) (2021) 7381.

\bibitem{yaji2016topology}
K.~Yaji, T.~Yamada, M.~Yoshino, T.~Matsumoto, K.~Izui, S.~Nishiwaki, Topology
  optimization in thermal-fluid flow using the lattice {Boltzmann} method,
  Journal of Computational Physics 307 (2016) 355--377.

\bibitem{xie2021topology}
S.~Xie, K.~Yaji, T.~Takahashi, H.~Isakari, M.~Yoshino, T.~Matsumoto, Topology
  optimization for incompressible viscous fluid flow using the lattice kinetic
  scheme, Computers \& Mathematics with Applications 97 (2021) 251--266.

\bibitem{suarez2018topology}
M.~A. Suarez, J.~S. Romero, I.~F. Menezes, {Topology Optimization for fluid
  flow problems using the Virtual Element Method}, Mec{\'a}nica Computacional
  36~(45) (2018) 2037--2046.

\bibitem{padhy2023fluto}
R.~K. Padhy, A.~Chandrasekhar, K.~Suresh, {FluTO: Graded} multi-scale topology
  optimization of large contact area fluid-flow devices using neural networks,
  Engineering with Computers (2023) 1--17.

\end{thebibliography}

\end{document}